\documentclass[10pt]{amsart}   
\usepackage{amssymb,amsmath,amscd,amsfonts,amsthm,mathrsfs, curves}
\usepackage[active]{srcltx} 
\usepackage[all]{xy}  
\usepackage{color}

\def\thesection{\arabic{section}} 
\renewcommand{\theequation}{\arabic{section}.\arabic{equation}} 
\newtheorem{theorem}{Theorem}[section] 
\newtheorem{lemma}[theorem]{Lemma} 
\newtheorem{proposition}[theorem]{Proposition} 
\newtheorem{corollary}[theorem]{Corollary} 
\newtheorem{assumption}[theorem]{Assumption} 
\theoremstyle{definition}  
\newtheorem{definition}[theorem]{Definition} 
 
\theoremstyle{remark}  
\newtheorem{rem}[theorem]{Remark}  
  
\newcommand{\Gr}{\mathscr{G}}
\newcommand{\Hr}{\mathscr{H}}

\newcommand{\Sr}{\mathscr{S}}

\newcommand{\Bc}{\mathcal{B}}  
\newcommand{\Cc}{\mathcal{C}}  
\newcommand{\Dc}{\mathcal{D}} 
\newcommand{\Fc}{\mathcal{F}}  
\newcommand{\Hc}{\mathcal{H}}  
\newcommand{\Ic}{\mathcal{I}}  
\newcommand{\Jc}{\mathcal{J}}

\newcommand{\Sc}{\mathcal{S}}

\newcommand{\Wc}{\mathcal{W}}  
\newcommand{\ad}{\mathrm{ad}}

\newcommand{\re}{\mathrm{Re}}   
\newcommand{\im}{\mathrm{Im}}

\newcommand{\supp}{\mathrm{supp}}  
\newcommand{\C}{\mathbb{C}}  
\newcommand{\R}{\mathbb{R}}  
\newcommand{\N}{\mathbb{N}}  
\newcommand{\Z}{\mathbb{Z}}

\def\build#1_#2^#3{\mathrel{\mathop{\kern 0pt#1}\limits_{#2}^{#3}}} 
\def\cchi{\raisebox{.45 ex}{$\chi$}}
 
\newcommand{\un}{1\hskip-0.6ex {\mathsf{ I} }} 

\newcommand{\op}{{\rm Op\, }}  
\newcommand{\ph}{{\rm x}^\ast}
\newcommand{\ssi}{\Longleftrightarrow}

\newcommand{\impl}{\Longrightarrow}
\newcommand{\donne}{\mapsto}
\newcommand{\dans}{\longrightarrow}
\newcommand{\rH}{{\rm H}}

\newcommand{\rL}{{\rm L}}                 

\begin{document}  
\title[Weighted Mourre theory]  
{Weighted Mourre's commutator theory, application to Schr\"odinger operators with 
oscillating potential}   
\author{Gol\'enia, Sylvain}   
\address{IMB, UMR 5251 du CNRS, Universit\'e de Bordeaux 1\\ 
351 cours de la Lib\'eration\\ 
F-33 405 Talence cedex, France} 
\email{Sylvain.Golenia@math.u-bordeaux1.fr}   
 \author{Jecko, Thierry}   
\address{AGM, UMR 8088 du CNRS, Universit\'e de Cergy Pontoise, Site de Saint Martin,  
2, rue Adolphe Chauvin, F-95000 Cergy-Pontoise, France}   
\email{thierry.jecko@u-cergy.fr}   
\date{\today}   
\begin{abstract}   
We present a variant of Mourre's commutator theory. We apply it to prove the limiting 
absorption principle for Schr\"odinger operators with a perturbed Wigner-Von Neumann potential 
at suitable energies. To our knowledge, this result is new since we allow 
a long range perturbation of the Wigner-Von Neumann
potential. Furthermore, we can show that  
the usual Mourre theory, based on differential inequalities and on the 
generator of dilation, cannot apply to the mentioned Schr\"odinger operators. 
\end{abstract}   
\maketitle   
\tableofcontents 

\section{Introduction.} \label{s:intro}
\setcounter{equation}{0}

Since its introduction in 1980 (cf., \cite{m}), many papers have shown
the power of Mourre's  commutator theory to study the point and
continuous spectra of a quite wide class of self-adjoint
operators. Among others, we refer to
\cite{bfs,bchm,cgh,dj,fh,ggm1,ggo,hus,jmp,s} and to the  book
\cite{abg}. One can also find parameter dependent versions of the
theory (a semi-classical one for  instance) in \cite{rt,w,wz}. Recently
it has been extended to (non self-adjoint) dissipative  operators
(cf., \cite{bg,roy}). \\ 
In \cite{gj}, we introduced a new approach of Mourre's commutator
theory, which is strongly inspired by results in semi-classical
analysis (cf., \cite{b,cj,jec1,jec2}). In \cite{ge}, C.\ G\'erard
showed that  it can be followed using traditional ``energy
estimates''. This approach furnishes an alternative way to develop the
original Mourre Theory and do not use 
differential  inequalities. \\
The aim of the present paper is to present a new theory, which 
shares common goals with Mourre's commutator theory but relies
on different assumptions.  
It is inspired by the approach in \cite{gj} and in \cite{ge}. It is
actually new  since we can produce an example for which it applies while the strongest versions 
of Mourre's commutator theory (cf., \cite{abg,s}) with (variants
of) the generator of dilation  as conjugate operator cannot be
applied to it. \\ 
Our example is a perturbation of a Schr\"odinger operator with a Wigner-Von Neumann potential. 
Furthermore we can allow a long range perturbation which is not covered by 
previous results in \cite{dmr,ret1,ret2}. A similar situation is considered in \cite{mu} 
but at different energies. 

Let us now briefly recall Mourre's commutator theory and present our results. We need some 
notation and basic notions (see 
Section~\ref{s:base} for details). We consider two self-adjoint 
(unbounded) operators $H$ and $A$ acting in some complex Hilbert space
$\Hr$. Let $\|\cdot\|$ denote the norm of bounded operators on $\Hr$. 
With the help of $A$, we study spectral properties of $H$, the spectrum 
$\sigma (H)$ of which is included in $\R$. 
Let $\Ic, \Jc$ be open intervals of $\R$. Given $k\in \N$, 
we say that $H\in \Cc_\Jc ^k(A)$ if for all $\cchi \in \Cc^\infty_c(\R)$ with support in $\Jc$, 
for all $f\in \Hr$, the map $\R\ni t\mapsto e^{itA}\cchi (H)e^{-itA}f\in \Hr$ has the usual 
$\Cc^k$ regularity. Denote by $E_\Ic(H)$ the spectral measure of $H$ above $\Ic$. We say that 
the \emph{Mourre estimate} holds true for $H$ on $\Ic$ if there exist $c>0$ and a compact 
operator $K$ such that   
\begin{eqnarray}\label{eq:mourre} 
E_\Ic (H)[H, iA]E_\Ic(H)\geq E_\Ic(H)\, (c\, +\, K)\, E_\Ic(H), 
\end{eqnarray}  
in the form sense on $\left(\Dc(A)\cap \Dc(H)\right) \times
\left(\Dc(A)\cap \Dc(H)\right)$. In general, the l.h.s.\ of  
\eqref{eq:mourre} does not extend, as a form, on $\Hr\times \Hr$ but
it is the case if $H\in \Cc_\Jc ^1(A)$  and $\Ic\subset\Jc$ (cf.,
\cite{s,gj}). We say that the  \emph{strict Mourre estimate} holds
true if the Mourre estimate \eqref{eq:mourre} holds  true with
$K=0$. In the first case (resp.\  the second case), it turns out that
the point  spectrum of $H$ is finite (resp.\  empty) in compact
sub-intervals $\Ic '$ of $\Ic$ if  $H\in \Cc_\Jc ^k(A)$ and
$\Ic\subset\Jc$.  The main aim of Mourre's commutator theory is to
show, when the strict Mourre estimate  holds true for $H$ on $\Ic$,
the following \emph{limiting absorption principle} (LAP)  on compact
sub-intervals $\Ic '$ of $\Ic$. Given such a $\Ic '$ and $s>1/2$, we
say that  the LAP, respectively to the triplet $(\Ic ',s,A)$, holds
true for $H$ if   
\begin{equation}\label{eq:lap}  
\sup_{{\rm Re}z\in \Ic ', {\rm Im}z\neq 0}\|\langle A\rangle^{-s}(H-z)^{-1}
\langle A\rangle^{-s}\|<\infty ,     
\end{equation}  
where $\langle t\rangle = (1+|t|^2)^{1/2}$. In that case, 
it turns out that the spectrum of $H$ is purely absolutely continuous
in $\Ic '$ (cf., Theorem XIII.20 in \cite{rs4}). 
Notice that \eqref{eq:lap} holds true for $s=0$ 
if and only if $\Ic '\cap \sigma (H)=\emptyset$. \\
In \cite{abg,s}, such LAPs are derived under a slightly stronger regularity assumption than 
$H\in \Cc_\Jc ^1(A)$ with $\Ic\subset\Jc$. Actually, stronger results are proved. In 
particular, in the norm topology of bounded operators, one can defined the boundary values 
of the resolvent: 
\begin{equation}\label{bord}
\Ic '\ni\lambda \donne \lim _{\varepsilon \to 0^\pm}\langle A\rangle^{-s}(H-\lambda -i\varepsilon)^{-1}
\langle A\rangle^{-s}
\end{equation}
and show some H\"older continuity for them. \\
Implicitly in \cite{gj} and explicitly in \cite{ge}, one 
can derive, using $H\in \Cc_\Jc ^2(A)$ with $\Ic\subset\Jc$, the LAP \eqref{eq:lap} on 
compact sub-intervals $\Ic '$ of $\Ic$ from the Mourre estimate
\eqref{eq:mourre} with $K=0$ via a \emph{strict, weighted Mourre estimate}: 
\begin{equation}\label{eq:resc-mourre} 
E_\Ic (H)[H, i\varphi (A)]E_\Ic (H)\ \geq \ 
c_1 E_\Ic (H)\langle A\rangle^{-1-\varepsilon}E_\Ic (H), 
\end{equation} 
where $\varepsilon =2s-1>0$ and $\varphi$ is some appropriate non-negative, bounded, smooth function 
on $\R$. Note that the l.h.s.\ of \eqref{eq:resc-mourre} is a well defined form on 
$\Hr\times \Hr$. It seems that the use of such kind of inequality to derive resolvent estimates 
appears in \cite{jec1} for the first time. \\
Our new idea is to take the strict, weighted Mourre estimate \eqref{eq:resc-mourre} as starting point, instead of the 
strict Mourre estimate. This costs actually less regularity of $H$ w.r.t.\ $A$. Precisely, we show 
\begin{theorem}\label{th:mourre} 
Let $\Ic$ be a bounded, open interval of $\R$ and assume that $H\in \Cc_{\Ic}^1(A)$. 
Assume that, for some $\varepsilon_0 >0$, for any $\varepsilon\in (0;\varepsilon_0]$, there exists some 
real borelian bounded function $\varphi$ such that the strict, weighted Mourre estimate, 
i.e. \eqref{eq:resc-mourre}, holds true. Then, for any $s>1/2$ and for any closed sub-interval 
$\Ic'$ of $\Ic$, the LAP \eqref{eq:lap} for $H$ respectively to $(\Ic', s, A)$ holds true. 
\end{theorem} 
\begin{rem}\label{r:1/2<s<1}
Notice that the LAP \eqref{eq:lap} for $H$ respectively to 
$(\Ic', s, A)$ implies the LAP \eqref{eq:lap} for $H$ respectively to 
$(\Ic', s', A)$, for any $s'\geq s$. Therefore, it is enough to prove
Theorem \ref{th:mourre} for $s$ close to $1/2$. 
\end{rem}
\begin{rem}\label{r:controle-borne}
Using G\'erard's energy method in \cite{ge}, we can upper bound the size of the l.h.s.\ of 
\eqref{eq:lap} in terms of the constant $c_1$ appearing in \eqref{eq:resc-mourre}. See 
Corollary~\ref{co:mourre-ameliore}. 
\end{rem}
Actually Theorem~\ref{th:mourre} will follow from the more general result obtained in 
Theorem~\ref{th:mourre-reech-gene-proj}. The new theory that we present here and that 
we call ``weighted Mourre theory'' is essentially a part of the variant of the Mourre 
theory in \cite{ge,gj}. As such, it is simpler than the usual Mourre theory (it does 
not use differential inequalities). However, we do not know if such approach gives 
continuity results on the boundary values of the resolvent 
\eqref{bord}. We shall give two (almost equivalent) ways to view the new theory 
(cf., Subsections~\ref{ss:special} and~\ref{ss:energy}). 

As announced above, we want to derive the LAP \eqref{eq:lap} (for some $A$) on
carefully chosen intervals $\Ic '$ for a certain class of Schr\"odinger
operators. Let $d\in\N^\ast$ and let $H_0$ be the self-adjoint
realization of the Laplacian $-\Delta _x$ in $\rL ^2(\R ^d_x)$. 
Given $q\in \R^\ast$ and $k>0$, the function $W:\R^d\dans\R$ defined by 
$W(x)=q(\sin k|x|)/|x|$ is called the Wigner-Von Neumann potential. We consider another 
real valued function $V$ satisfying some long range condition (see Section~\ref{s:wigner}
for details) such that the operator $H_1:=H_0+W+V$ is self-adjoint on
the domain of $H_0$. This is the Schr\"odinger operator with
a perturbed Wigner-Von Neumann potential that we consider. It is well known that its essential spectrum 
is $[0;+\infty [$. Now we look for an interval $\Ic '\subset ]0;+\infty [$ 
on which we can get the LAP \eqref{eq:lap}. As operator $A$, it is natural to choose 
the generator of dilation $A_1$, the self-adjoint realization of $(x\cdot \nabla _x+\nabla
_x\cdot x)/(2i)$ in $\rL^2(\R^d_x)$. Indeed, when $W$ is absent, such LAPs have been derived. 
As mentioned above, the pure point spectrum $\sigma _{pp}(H_1)$ of $H_1$ has to be empty 
in $\Ic '$.\\ 
There are many papers on the absence of positive eigenvalue for Schr\"odinger operators: see 
\cite{k,si,a,fhhh2,fh,ij,rs4,cfks}. They do not apply to the present situation 
because of the behaviour of the Wigner-Von Neumann $W$. One can even show that $k^2/4$ 
is actually an eigenvalue of $H_1$ for a well chosen, radial, short range potential $V$ (cf., 
\cite{rs4} p.\ 223 and \cite{bd}). \\
In dimension $d=1$, the eigenvalue at $k^2/4$ is preserved under
suitable perturbation  
(see \cite{chm}). Furthermore it is proved in 
\cite{fh,fhhh1} that, if $|q|<k$, the usual Mourre estimate 
\eqref{eq:mourre} holds true on compact intervals 
$\Ic \subset ]0;+\infty[$ and there is no eigenvalue in $]0;+\infty[$,
and otherwise that,  
 on compact $\Ic \subset ]0;+\infty[\setminus\{k^2/4\}$, no eigenvalue
 is present and  
the usual Mourre estimate \eqref{eq:mourre} holds true. Actually if 
$k^2/4$ is an eigenvalue of $H_1$ then the usual Mourre estimates
cannot hold true on  a compact neighbourhood of $k^2/4$, with the
generator of dilation as conjugate operator. This follows
from the arguments of the proof of  Corollary 2.6 in \cite{fh}. Thus
the eigenvalue $k^2/4$ is a threshold. \\ 
We focus on compact intervals $\Ic$ satisfying $\{0, k^2/4\}\cap  \Ic=\emptyset$ 
and, when $d>1$, $\Ic \subset ]0;k^2/4 [$.
Using pseudodifferential calculus and recycling arguments from
\cite{fh}, we prove the usual Mourre estimate \eqref{eq:mourre} on
such $\Ic$, the operator $A$ being $A_1$, yielding the finiteness of the pure point spectrum  
$\sigma _{pp}(H_1)$ in $\Ic$.  Then, in Theorem \ref{th:tal-Wigner},  we derive a strict,
weighted Mourre estimate  \eqref{eq:resc-mourre} and show that
Theorem~\ref{th:mourre} applies, leading to the LAP \eqref{eq:lap} for $H_1$.  
For short range perturbation $V$, we partially recover
results from \cite{dmr,ret1,ret2} but, in contrast to these papers, we
are able to treat a long range perturbation $V$. We mention that in
\cite{mu}, for high  enough energies, one proves a LAP for long-range
perturbations of a larger class of oscillating potentials. In this
situation, the Laplacian is a priori ``stronger'' than the potential,
in contrast to  
the present case. \\
As already mentioned, the LAP  for $H_1$ implies 
the absence of singular spectrum over $\Ic$. This result seems to be
new, even in dimension $1$. Concerning this question for potentials,
which are ``decaying at most like $1/x$ at infinity'', we refer
to \cite{K,r} and  references therein. 
\\
Finally we show that $H_1$ does not have the required
regularity w.r.t.\ (variants of) $A_1$ to
apply the usual Mourre theory from  
\cite{abg,ggm1,s}. For the same reason, the derivation of 
the strict, weighted Mourre estimate \eqref{eq:resc-mourre} for $H_1$ from the 
corresponding strict Mourre estimate, i.e. \eqref{eq:mourre} with $K=0$, along the lines 
in \cite{ge}, is not allowed. If one removes the oscillating potential
$W$ and keeps our  assumptions on $V=V_{\rm lr}+V_{\rm sr}$ (see
  Assumption~\ref{a:long-range-3}), the situation is
    well-known. However, to obtain a LAP with the traditionnal Mourre theory,
one needs some strong and quite involved versions as in
\cite{abg,s}. In particular, the results in  
\cite{m,ge} do not apply in this case, while our arguments here gives
the LAP \eqref{eq:lap} on any  compact interval $\Ic \subset
]0;+\infty[$ (cf., Remarks~\ref{r:above-threshold}  
and~\ref{r:TAL-d=1->seuil} below). \\
We did not optimize our study of Schr\"odinger operators with oscillating potential. We 
believe that we can handle more general perturbations. Because of a
difficulty explained  in Remark~\ref{r:above-threshold}, we did not
consider intervals 
$\Ic$ above $k^2/4$ for $d>1$.  
However we believe that a variant of the present theory is applicable in this case. 
We think that a general study of long range perturbations of the Schr\"odinger operator 
with Wigner-Von Neumann potential is interesting in itself and hope to develop it in a forthcoming 
paper. \\
The paper is organized as follows. 
In Section~\ref{s:base}, we introduce some notation and 
basic but important notions. 
In Section~\ref{s:nouv-strategie}, we show a stronger version of 
Theorem~\ref{th:mourre}, namely
Theorem~\ref{th:mourre-reech-gene-proj}. 
In Section~\ref{s:wigner}, we study Schr\"odinger operators with
perturbed Wigner-Von Neumann potentials. In Subsection~\ref{ss:wigner-mourre}, we
derive usual Mourre estimates below the ``threshold'' $k^2/4$. In 
Subsection~\ref{ss:wigner-mourre-reech}, we essentially apply 
Theorem~\ref{th:mourre} to Schr\"odinger operators.
In Section~\ref{s:usual}, we prove that they cannot be treated by the usual Mourre theory 
in \cite{abg,s}. In Appendix~\ref{app:s:oscillation}, we prove a key pseudodifferential 
result to control the behaviour of the Wigner-Von Neumann potential
(extending a result  by \cite{fh} in dimension one). In
Appendix~\ref{app:s:calcul-fonct-pseudo}, we review functional
calculus for pseudodifferential operators (cf., \cite{bo}). In
Appendix~\ref{app:s:interpolation}, we establish the boundedness of
some operator  using interpolation. Finally, in
Appendix~\ref{app:s:dim1}, we present, in dimension one,  a simpler
proof of Lemma~\ref{l:non-compact}, this lemma being used to show that
the regularity assumption of the usual Mourre theory is not satisfied
by the Schr\"odinger operators studied here.

{\bf Acknowledgement:} The authors thank Jean-Michel Bony,
Vladimir Georgescu, Ira Herbst, Andreas Knauf, Jacob Schach M\o ller, Nicolas Lerner, 
Karel Pravda-Starov, and Erik Stibsted for fruitful discussions. \\ 
The authors apologize to Jacob Schach M\o ller for not citing his paper \cite{mo}
in their previous work \cite{gj} when they proved Proposition 2.6 and 2.7 below.
They did not realize the presence of the corresponding result in \cite{mo}.

\section{Basic notions and notation.}\label{s:base}
\setcounter{equation}{0}

In this section, we introduce some notation and recall known results. For details, we refer 
to \cite{abg,dg,gj,s} on regularity and to \cite{h3,bo,bo2,bc,l} on pseudodifferential 
calculus.

\subsection{Regularity.} 
\label{ss:regularity}

For an interval $\Ic$ of $\R$, we denote by $\overline{\Ic}$
(resp.\  $\mathring{\Ic}$) its closure (resp.\  its interior). 
The scalar product $\langle \cdot ,\cdot \rangle$ in $\Hr$ is right linear and $\|\cdot\|$ 
denotes the corresponding norm and also the norm in $\Bc (\Hr )$, the space of bounded 
operators on $\Hr$. Let $A$ be a self-adjoint operator. Let $T$ be a closed operator. The 
form $[T,A]$ is defined on $(\Dc(A)\cap\Dc(T))\times(\Dc(A)\cap\Dc(T))$ by 
\begin{equation}\label{eq:sens-forme}
 \langle f\, ,\, [T,A]g\rangle \ :=\ \langle T^\ast f\, ,\, Ag\rangle\, -\, 
\langle Af\, ,\, Tg\rangle \, .
\end{equation}
If $T$ is a bounded operator on $\Hr$ and $k\in \N$, we say that $T\in \Cc^k(A)$ 
if, for all $f\in \Hr$, the map $\R\ni t\mapsto e^{itA}Te^{-itA}f\in \Hr$ has the usual 
$\Cc^k$ regularity. The following characterization is available. 
\begin{proposition}\label{p:caract-C1} (\cite[p.\ 250]{abg}). 
Let $T\in \Bc(\Hr)$. Are equivalent: 
\begin{enumerate}
 \item $T\in \Cc^1(A)$.  
 \item The form $[T,A]$ defined on $\Dc(A)\times \Dc(A)$ extends to 
a bounded form on $\Hr\times\Hr$ associated to a bounded operator denoted by  $\ad_A^1(T):=[T,A]_\circ$. 
 \item $T$ preserves $\Dc (A)$ and 
the operator $TA-AT$, defined on $\Dc(A)$, extends to a bounded operator on $\Hr$.
\end{enumerate}
\end{proposition} 
It follows that $T\in \Cc^k(A)$ if and only if the iterated
commutators $\ad_A^p(T):= 
[\ad_A^{p-1}(T),A]_\circ$ are bounded for $p\leq k$. In particular,
for $T\in \Cc^1(A)$,  $T\in \Cc^2(A)$ if and only if $[T,A]_\circ\in
\Cc^1(A)$. \\ 
Let $H$ be a self-adjoint operator and $\Ic$ be an open interval. As
in the Introduction  (Section~\ref{s:intro}), we say that $H$ is {\em
  locally of class $\Cc^k(A)$ on $\Ic$},  we write $H\in
\Cc^k_\Ic(A)$, if, for all $\varphi\in\Cc_c^\infty(\Ic)$,  
$\varphi(H)\in \Cc^k(A)$. \\
It turns out that $T\in \Cc^k(A)$ if and only if, for a $z$ outside
$\sigma (T)$, the spectrum of $T$, $(T-z)^{-1}\in \Cc^k(A)$. It  is
natural to say that $H\in \Cc^k(A)$ if $(H-z)^{-1}\in \Cc^k(A)$ for
some  $z\not\in\sigma (H)$. In that case, $(H-z)^{-1}\in \Cc^k(A)$,
for all $z\not\in\R$. This regularity is stronger than the local one
as asserted in the following  
\begin{proposition}\label{p:global-impl-local} (\cite[p.\ 244]{abg})
If $H\in\Cc^k(A)$ then $H\in \Cc^k_\Ic(A)$ for all open interval $\Ic$ of $\R$. 
\end{proposition}
Next we recall Proposition~2.1 in \cite{gj} which gives a sufficient condition to get 
the $\Cc^1(A)$ regularity for finite range operators. 
\begin{proposition}\label{p:rang-fini} (\cite{gj})
If $f,g\in \Dc (A)$, then the rank one operator 
$|f\rangle \langle g|:h\mapsto \langle g,h\rangle f$ is in $\Cc^1(A)$. 
\end{proposition}
For $\rho\in\R$, let $\Sc^\rho$ be the class of functions $\varphi\in\Cc^\infty(\R)$ 
such that  
\begin{eqnarray}\label{eq:def-S-rho} 
\forall k\in\N, \quad C_k(\varphi) :=\sup _{t\in\R}\, \langle t\rangle^{-\rho+k}|
\varphi^{(k)}(t)|<\infty .     
\end{eqnarray} 
Here $\varphi^{(k)}$ denotes the $k$th derivative of $\varphi$. 
Equipped with the semi-norms defined by (\ref{eq:def-S-rho}), $\Sc^\rho$ is
a Fr\'echet space. We recall the following result from \cite{dg} on almost analytic extension. 
\begin{proposition}\label{p:dg}(\cite{dg})
Let $\varphi\in\Sc^\rho$ with $\rho\in\R$. There is a smooth function 
$\varphi^\C:\C \rightarrow \C$, called an 
\emph{almost analytic extension} of $\varphi$, such that, for all   
$l\in \N$,
\begin{align} 
\label{eq:dg1} \varphi^\C|_{\R}=\varphi,\quad &\big|\partial_{\overline{z}}\varphi^\C(z)
\big|\leq c_1 \langle \re(z)\rangle^{\rho-1 -l} |\im(z)|^l\, ,\\
\label{eq:dg2} 
& \supp\, \varphi^\C \subset\{x+iy; |y|\leq c_2 \langle 
  x\rangle\},\\
\label{eq:dg3} & \varphi^\C(x+iy)= 0, \mbox{ if } 
  x\not\in\supp\,\varphi,  
\end{align} 
for constants $c_1$, $c_2$ depending on the semi-norms 
\eqref{eq:def-S-rho} of $\varphi$ in $\Sc^\rho$.   
\end{proposition}  
%


Next we recall Helffer-Sj\"{o}strand's functional calculus
(cf., \cite{hs,dg}). For $\rho< 0$, $k\in\N$, and $\varphi\in \Sc^{\rho}$, the bounded 
operators $\varphi ^{(k)}(A)$ can be recovered by 
\begin{eqnarray}\label{eq:int} 
\varphi ^{(k)}(A) = \frac{i(k!)}{2\pi}\int_\C\partial_{\overline{z}}\varphi^\C(z)(z-A)^{-1-k}
dz\wedge d\overline{z}, 
\end{eqnarray} 
where the integral exists in the norm topology, by \eqref{eq:dg1} 
with $l=1$. For $\rho\geq 0$, we rely on the following approximation:
\begin{proposition}\label{p:ext-hs}(\cite{gj})
Let $\rho \geq 0$ and $\varphi\in \Sc^{\rho}$. Let $\cchi \in\Cc
_c^\infty(\R)$ with $\cchi =1$ near $0$ and $0\leq\cchi \leq 1$, and,
for $R>0$, let $\cchi _R(t)=\cchi (t/R)$. For $f\in
\Dc(\langle A\rangle^{\rho})$, there exists 
\begin{eqnarray}\label{eq:int-cv} 
\varphi ^{(k)}(A)f =\lim_{R\to +\infty}\frac{i}{2\pi}\int_\C\partial_{\overline{z}}
(\varphi \cchi_R)^\C(z)
(z-A)^{-1-k}f\, dz\wedge d\overline{z}.
\end{eqnarray}  
The r.h.s.\ converges for the norm in $\Hr$. It is independent of the
choise of $\cchi$.
\end{proposition}  

Notice that, for some $c>0$ and $s\in [0;1]$, there exists
some $C>0$ such that, for all $z=x+iy\in\{a+ib\mid 0<|b|\leq c\langle a\rangle \}$ 
(like in \eqref{eq:dg2}), 
\begin{eqnarray}\label{eq:majoA} 
\big\| \langle A\rangle^s (A-z)^{-1}\big\|\leq C \langle x \rangle^{s}\cdot |y|^{-1}. 
\end{eqnarray}   

Observing that the self-adjointness assumption on $B$ is useless, we pick from \cite{gj,mo} the 
following result in two parts. 
\begin{proposition}\label{p:regu}(\cite{gj,mo})
Let $k\in \N^\ast$, $\rho< k$, $\varphi\in \Sc^{\rho}$, and $B$ be a bounded
operator in $\Cc^k(A)$. As forms on 
$\Dc(\langle A\rangle^{k-1})\times \Dc(\langle A\rangle^{k-1})$,  
\begin{align}\label{eq:egalite}
[\varphi(A), B] &= \sum_{j=1}^{k-1} \frac{1}{j!} 
\varphi^{(j)}(A)\ad_A^j(B)\\ 
\label{eq:reste22}    
&\quad + \,
\frac{i}{2\pi}\int_\C\partial_{\overline{z}}\varphi^\C(z)(z-A)^{-k}
\ad_A^k(B) 
(z-A)^{-1} dz\wedge d\overline{z}. 
\end{align} 
In particular, if $\rho\leq 1$, then $B\in\Cc^1(\varphi(A))$. 
\end{proposition} 

The rest of the previous expansion is estimated in 
\begin{proposition}\label{p:est3} (\cite{gj,mo})
Let $B\in\Cc^k(A)$ bounded. Let $\varphi\in\Sc^\rho$, with $\rho< k$. Let 
$I_k(\varphi)$ be the rest of the development of order $k$ \eqref{eq:egalite} of $[\varphi(A), B]$, 
namely \eqref{eq:reste22}. Let 
$s, s'\geq 0$ such that $s'<1$, $s<k$, and $\rho+s+s'<k$. Then, for $\varphi$ staying in a 
bounded subset of $\Sc^\rho$, $\langle A \rangle^{s} 
I_k(\varphi)\langle A \rangle^{s'}$ is bounded and there exists a $A$ and $\varphi$
independent constant $C>0$ such that $\|\langle A \rangle^{s} 
I_k(\varphi)\langle A \rangle^{s'}\|\leq C\|\ad_A^k(B)\|$. 
\end{proposition}
We refer to \cite{bg} for some generalization of Propositions
\ref{p:regu} and \ref{p:est3} to the case where $B$ is unbounded and
$[A, B]_\circ$ is bounded.

\subsection{Pseudodifferential calculus.} 
\label{ss:psi-do}

In this subsection, we briefly review some basic facts about
pseudodifferential calculus that  
we need in the treatment of Schr\"odinger operators. 
We refer to \cite{h3}[Chapters 18.1, 18.4, 18.5, and 18.6] for a
traditional study  
of the subject but also to \cite{bo,bo2,bc,l} for a modern and powerful version. Other results are 
presented in Appendix~\ref{app:s:oscillation} 
and~\ref{app:s:calcul-fonct-pseudo}. 

Denote by $\Sc (M)$ the Schwartz space on the space $M$ and by $\Fc$ the Fourier transform 
on $\R^d$ given by  
\[\Fc u (\xi) :=  (2\pi )^{-d}\int_{\R^{d}}e^{-ix\cdot\xi}u(x)\, dx\, ,\]
for $\xi\in\R^{d}$ and $u\in\Sc (\R^d)$. For test functions $u, v\in\Sc (\R^d)$, 
let $\Omega (u,v)$ and $\Omega '(u,v)$ be the functions in $\Sc
(\R^{2d})$ defined by  
\begin{align*}
\Omega (u,v)(x,\xi )&:=\overline{v}(x)\Fc u(\xi )e^{ix\cdot\xi}\, ,\\
\Omega '(u,v)(x,\xi )&:=(2\pi
)^{-d}\int_{\R^{d}}u(x-y/2)\overline{v}(x+y/2)e^{-iy\cdot\xi}\, dy\, , 
\end{align*}
respectively. Given a distribution $b\in \Sc '(T^\ast\R^d)$, the formal quantities 
\begin{eqnarray*}
(2\pi )^{-d}\int_{\R^{3d}}e^{i(x-y)\cdot\xi}b(x,\xi)v(x)u(y)\, dxdyd\xi\, , \\
(2\pi )^{-d}\int_{\R^{3d}}e^{i(x-y)\cdot\xi}b((x+y)/2,\xi)u(x)u(y)\, dxdyd\xi\,
\end{eqnarray*}
are defined by the duality brackets $\langle b,\Omega (u,v)\rangle$ and 
$\langle b,\Omega '(u,v)\rangle$, respectively. They define continuous operators 
from $\Sc (\R^d)$ to $\Sc '(\R^d)$ that we denote by $\op b(x,D_x)$ 
and $b^w(x,D_x)$ respectively. Sometimes we simply write $\op b$ and
$b^w$, respectively.\\
Choosing on the phase space $T^\ast\R^d$ a metric $g$ and a weight function $m$ with 
appropriate properties (cf., admissible metric and weight in \cite{l}), let 
$S(m,g)$ be the space of smooth functions on $T^\ast\R^d$ such that, for all $k\in \N$, 
there exists $c_k>0$ so that, for all $\ph =(x,\xi )\in T^\ast\R^d$, all $(t_1,\cdots ,t_k)
\in (T^\ast\R^d)^k$, 
\begin{equation}\label{eq:bound-symbol}
|a^{(k)}(\ph )\cdot (t_1,\cdots ,t_k)|\leq c_k m(x^*)g_{x^*}(t_1)^{1/2}\cdots 
g_{x^*}(t_k)^{1/2}\, .
\end{equation}
Here, $a^{(k)}$ denotes the $k$-th derivative of $a$. 
We equip the space $S(m,g)$ with the semi-norms $\|\cdot\|_{\ell ,S(m,g)}$
defined by $\max_{0\leq k\leq\ell}c_k$, where the $c_k$ are 
the best constants in \eqref{eq:bound-symbol}. $S(m,g)$ is a Fr\'echet space. The space 
of operators $\op b(x,D_x)$ (resp.\  $b^w(x,D_x)$) when $b\in S(m,g)$ has nice properties 
(cf., \cite{h3,l}). Defining $\ph =(x,\xi )\in T^\ast\R^d$, we stick here to the following metrics 
\begin{equation}\label{eq:metric}
 g_{\ph}:=\ \frac{dx^2}{\langle x\rangle ^2}+\frac{d\xi ^2}{\langle \xi\rangle ^2}
\hspace{.5cm}\mbox{and}\hspace{.5cm} (g_0)_{\ph}:=\ dx^2+\frac{d\xi ^2}{\langle \xi\rangle ^2},
\end{equation}
and to weights of the form, for $p, q \in \R$, 
\begin{align}\label{eq:weight}
m(\ph ):=\langle x\rangle ^p\langle \xi\rangle ^q.
\end{align}
The gain of the calculus associated to each metric in \eqref{eq:metric} 
is given respectively by 
\begin{align}\label{eq:gain}
h(\ph ):=\langle x\rangle ^{-1}\langle \xi\rangle ^{-1} \mbox{ and } 
h_0(\ph )=\langle \xi\rangle ^{-1}. 
\end{align}
We note that $S(m,g)\subset S(m,g_0)$ with continuous injection.
Take weights $m_1$, $m_2$ as in \eqref{eq:weight}, let $\tilde g$ be $g$ or $g_0$, and denote 
by $\tilde h$ the gain of $\tilde g$. For any $a\in S(m_1,\tilde g)$ 
and $b\in S(m_2,\tilde g)$, there are a symbol $a\# _rb\in S(m_1m_2,\tilde g)$ and a symbol 
$a\# b\in S(m_1m_2,\tilde g)$ such that $\op a\op b=\op (a\# _rb)$ and 
$a^wb^w=(a\# b)^w$. The maps $(a,b)\donne a\# _rb$ and $(a,b)\donne a\# b$ 
are continuous and so are also $(a,b)\donne a\# _rb-ab\in
S(m_1m_2\tilde h,\tilde g)$ and  
$(a,b)\donne a\# b-ab\in S(m_1m_2\tilde h,\tilde g)$. If $a\in
S(m_1,\tilde g)$, there exists  
$c\in S(m_1,\tilde g)$ such that $a^w=\op c$. The maps $a\donne c$ and
$a\donne c-a 
\in S(m_1m_2\tilde h,\tilde g)$ are continuous. If $a\in S(1,\tilde
g)$, $a^w$ and  $\op a$ are bounded on $\rL^2(\R^d)$. For $a\in
S(1,\tilde g)$,  
\begin{equation}\label{eq:caract-pseudo-compact}
\op a\mbox{ is compact}\ssi a^w\mbox{ is compact}\ssi \lim _{|\ph|\to
  \infty}a(\ph)=0\, . 
\end{equation}

\section{Weighted Mourre theory.} \label{s:nouv-strategie}
\setcounter{equation}{0}

In this section, we present our new strategy to get the LAP
(\ref{eq:lap}). As in \cite{gj} (see also \cite{cgh}), we consider a more general version 
of the LAP, namely the LAP for the reduced resolvent (see \eqref{eq:tal-reduit} below). 
First we make use of a kind of weighted Weyl sequence introduced in \cite{gj}, that 
we call ``special sequence''. Then we present an adapted version of the method introduced in 
\cite{ge} and based on energy estimates. Both methods are quite close, the latter having the 
advantage to give an idea of the size of the l.h.s.\ of  (\ref{eq:lap}) (resp.\  \eqref{eq:tal-reduit}).

\subsection{Reduced resolvent.} 
\label{ss:reduced}

Let $P$ be the orthogonal projection onto the pure point spectral subspace of $H$ and $P^\perp =1-P$. 
For $s\geq 0$ and $\Ic '$ an interval of $\R$, we say that the {\em reduced} LAP, respectively to the triplet 
$(\Ic ',s,A)$, holds true for $H$ if  
\begin{equation}\label{eq:tal-reduit}  
 \sup_{{\rm Re}z\in \Ic ', {\rm Im}z\neq 0}\|\langle A \rangle^{-s}(H-z)^{-1} P^\perp
\langle A\rangle^{-s}\|<\infty .
\end{equation} 
Let $\Ic$ be an interval in $\R$ containing $\Ic '$ in its interior. 
Since $(H-z)^{-1}(1-E_\Ic(H))$ is uniformly bounded for $\re (z)\in \Ic '$ and $\im (z)\neq 0$, \eqref{eq:tal-reduit} is equivalent to the same estimate with $P^\perp$ replaced by $E_\Ic(H)P^\perp$. 
If no point spectrum is present in $\Ic$, then $(H-z)^{-1}E_\Ic(H)P^\perp =(H-z)^{-1}$
for ${\rm Re} z\in\Ic$ and (\ref{eq:tal-reduit}) is equivalent to the usual LAP (\ref{eq:lap}). \\
In \cite{cgh} and more recently in \cite{fms2}, it is shown that the
reduced LAP can be derived from the Mourre estimate
(\ref{eq:mourre}). In this case, it is well known that  the point
spectrum of $H$ is finite in $\Ic$ (but non empty in general, see
\cite{abg,m}).  In \cite{gj}, the reduced LAP is deduced from a
projected version of this Mourre estimate, namely 
\begin{equation}\label{eq:esti-mourre-proj}  
 P^\perp E_\Ic (H)[H, iA]  E_\Ic(H) P^\perp\geq c E_\Ic(H)P^\perp +
 P^\perp KP^\perp , 
\end{equation} 
for some compact operator $K$. In the proofs, one uses the
compactness of $K$ and the fact that the strong limit
\begin{equation}\label{eq:lim-fo}  
s-\lim_{\delta\to0} E_{]\lambda -\delta ;\lambda +\delta [}(H)P^\perp =0, 
\end{equation} 
to derive from \eqref{eq:esti-mourre-proj} a strict Mourre estimate (with
$K=0$) on all small enough intervals inside $\Ic$. Notice that the traditional 
theory (cf., \cite{abg,m}) performs the same derivation. So both methods rely on some strict Mourre 
estimate. Here, to get the reduced LAP \eqref{eq:tal-reduit} as shown in 
Theorem~\ref{th:mourre-reech-gene-proj} below, we also starts from a convenient strict estimate 
namely a strict, weighted, projected Mourre estimate like \eqref{eq:resc-mourre}. We discuss 
the possibility to derive it from a more general one in Subsection~\ref{ss:application}. Since 
we work with projected estimates, we need some regularity of $P^\perp$ w.r.t.\ $A$.

\subsection{Special sequences.} 
\label{ss:special}

We work in a larger framework. 
\begin{definition} \label{suite-speciale-abs}
Let $C$ be an injective, bounded, self-adjoint operator. Let $\Ic '$
be an interval of $\R$. 
\begin{enumerate}
\item 
A \emph{special sequence} $(f_n, z_n)_{n\in \N}$ for $H$ associated to
$(\Ic ',C)$ is a sequence $(f_n,z_n)_n\in (\Dc(H)\times \C)^\N$ such that, for
some $\eta \geq 0$, $\re(z_n)\in \Ic '$, $0\neq\im(z_n)\rightarrow 0$, 
$\|Cf_n \|\rightarrow \eta$, $P^\perp f_n=f_n$, $(H-z_n)f_n\in\Dc(C^{-1})$, 
and $\|C^{-1}(H-z_n)f_n \|\rightarrow 0$. The limit $\eta$ is called the 
\emph{mass} of the special sequence. 
\item The {\em reduced} LAP, respectively to $(\Ic',C)$, holds true for $H$ if 
\begin{equation}\label{eq:r-lap-C}  
\sup_{{\rm Re}z\in \Ic ', {\rm Im}z\neq 0}\|C(H-z)^{-1}P^\perp C\|<\infty .
\end{equation}  
\end{enumerate}
\end{definition}  
Notice that \eqref{eq:r-lap-C} for $C=\langle A\rangle ^{-s}$ with $s\in ]1/2;1[$ gives the 
LAP \eqref{eq:tal-reduit}, thanks to Remark~\ref{r:1/2<s<1}. 
\begin{proposition}\label{p:non-lap-localisation-suite-spe} 
Let $\Ic '$ be an interval of $\R$. Take an injective, bounded, self-adjoint 
operator $C$ such that, for some $\cchi$, a bounded, borelian function on
$\R$ with $\cchi =1$ near $\Ic '$, the operator $C\cchi (H)P^\perp C^{-1}$ extends 
to a bounded operator. Let $\theta$ be a borelian function on $\R$ such that
$\theta \cchi =\cchi$. Then the reduced LAP \eqref{eq:r-lap-C} holds
true if and only if, for all special sequence $(f_n, z_n)_n$ for $H$ associated to 
$(\Ic ', C)$ such that $\theta (H)f_n=f_n$ for all $n$, the corresponding mass is zero. 
\end{proposition} 
\proof Assume the LAP \eqref{eq:r-lap-C} true. Then, for any special 
sequence $(f_n, z_n)_n$ for $H$ associated to $(\Ic ', C)$, for all $n$, 
\[\|Cf_n\|\ \leq \ \|C(H-z_n)^{-1}P^\perp C\| \cdot \|C^{-1}(H-z_n)f_n \|,\]
yielding $\eta=0$. Now assume the LAP \eqref{eq:r-lap-C} false. Then there exists some complex 
sequence $(z_n)$ such that $\re z_n\in \Ic'$, $\im z_n\to 0$, and $\|C(H-z_n)^{-1}P^\perp C\|\to
\infty$. Since $((H-z_n)^{-1}(1-\cchi )(H))$ is uniformly bounded, 
we can find, for all $n$, nonzero $u_n\in \Hc$ and $0<\kappa _n\to 0$ such that 
\[\|C(H-z_n)^{-1}\cchi (H)P^\perp Cu_n\|\ =\ \|u_n\|/\kappa _n .\]
We set $f_n:=\kappa _n(H-z_n)^{-1}\cchi (H)P^\perp Cu_n/\|u_n\|$. Notice that $\theta
(H)f_n=P^\perp f_n=f_n$ and $\|Cf_n\|=1$. Since $C\cchi (H)P^\perp C^{-1}$ is bounded,
$\cchi (H)P^\perp $ preserves $\Dc (C^{-1})$, the image of $C$. Thus
$(H-z_n)f_n\in \Dc (C^{-1})$, for all $n$. We conclude by noticing that 
$\|C^{-1}(H-z_n)f_n\|\ \leq \ \kappa _n\cdot \|C^{-1}\cchi (H)P^\perp
C\|\ =\ o(1).$  \qed 
\begin{proposition}\label{p:virial-abs}  
Let $\Ic ', C$ be as in Proposition~\ref{p:non-lap-localisation-suite-spe}. 
Let $(f_n,z_n)_n$ be a special sequence for a self-adjoint operator
$H$ associated to  
$(\Ic ', C)$. For any \emph{bounded} operator $B$, such that
$CBC^{-1}$ extends to  a bounded operator, 
$\lim_{n \rightarrow \infty} \langle f_n, [H, B] f_n\rangle =0$. 
\end{proposition}  
\proof Since $(f_n,z_n)_n$ is a special sequence
and $CBC^{-1}$ is bounded, we obtain that $\langle (H-z_n)f_n\, ,\,
Bf_n\rangle =o(1)$ and $\langle (H-z_n)f_n\, ,\, f_n\rangle =o(1)$. 
Therefore, $2i {\rm Im}z_n \|f_n\|^2={\rm Im}\langle (H-z_n)f_n\, ,\,
f_n\rangle =o(1)$. Hence
\begin{align*}
\langle f_n\, ,\, [H , iB]f_n\rangle &=\langle
(H-\overline{z_n})f_n\, ,\, iBf_n\rangle \, -\,  \langle
B^\ast f_n\, ,\, i(H-z_n)f_n\rangle 
\\
&=-2i {\rm Im}z_n \cdot \langle
f_n\, ,\, iBf_n\rangle \, -\, 2{\rm Im}\langle
(H-z_n)f_n\, ,\, Bf_n\rangle =o(1). \hspace{.6cm}\qed
\end{align*}
\begin{theorem}\label{th:mourre-reech-gene-proj} 
Let $\Ic$ be an open interval and $\Ic'$ be a closed sub-interval of $\Ic$. Let $B,C$ 
be two bounded self-adjoint operators, $C$ being injective. Assume that, for some 
bounded, borelian function $\cchi$ on $\R$ with $\cchi =1$ on $\Ic '$
and $\supp\, \cchi\subset \Ic$, $C\cchi (H)P^{\perp}C^{-1}$ and $CBC^{-1}$ extend 
to bounded operators. Assume further that the following strict weighted projected Mourre estimate
\begin{equation}\label{eq:mourre-reech-theta-chi} 
P^\perp E_\Ic (H) [H , iB]E_\Ic (H)P^\perp\ \geq \ P^\perp E_\Ic (H)C^2E_\Ic (H)P^\perp 
\end{equation} 
is satisfied. Then the LAP \eqref{eq:r-lap-C} on $\Ic'$ holds true. 
\end{theorem} 
\proof Let $(f_n,z_n)_n$ be a special sequence for $H$ associated to 
$(\Ic ',C)$ such that $E_\Ic (H)f_n=f_n$ for all $n$. By 
Proposition~\ref{p:non-lap-localisation-suite-spe}, it suffices to
show that the mass $\eta$ of the special sequence is zero. Letting 
(\ref{eq:mourre-reech-theta-chi}) act on both sides on $f_n$, we infer that 
%
$\langle f_n\, ,\, [H , iB]f_n\rangle \ \geq \ 
\bigl\| Cf_n \bigr\|^2$.
%
Now Proposition~\ref{p:virial-abs} yields $\eta =0$. \qed

\proof[Proof of Theorem~\ref{th:mourre}] Thanks to Remark~\ref{r:1/2<s<1}, we may assume that 
$s\in ]1/2;1[$ with $\varepsilon :=2s-1\in ]0;\varepsilon_0]$. If, for $f\in\Dc (H)$ and 
for $E\in\Ic$, $Hf=Ef$ then, by \eqref{eq:resc-mourre}, 
$0\geq c_1\|\langle A\rangle ^{-(1+\varepsilon )/2}f\|^2$ and $f=0$. Thus 
$E_\Ic(H)P^\perp=E_\Ic(H)$ and \eqref{eq:resc-mourre} may be rewritten as 
\eqref{eq:mourre-reech-theta-chi} with $B=\varphi (A)$ and 
$C=\sqrt{c_1}\langle A\rangle^{-(1+\varepsilon)/2}=\sqrt{c_1}\langle A\rangle^{-s}$. 
Notice that the function of $A$, given by $CBC^{-1}$, extends to a bounded operator. Let 
$\cchi\in\Cc _c^\infty(\Ic)$ such that 
$\cchi =1$ on $\Ic '$. Since 
$\cchi (H)P^{\perp}=\cchi (H)E_\Ic(H)P^{\perp}=\cchi (H)E_\Ic(H)=\cchi (H)$ and 
$H\in \Cc_{\Ic}^1(A)$, $\cchi (H)P^{\perp}\in \Cc^1(A)$. By Proposition~\ref{p:regu}, 
$[\cchi (H)P^{\perp},\langle A\rangle^{s}]$ extends to a bounded operator. Thus, so does 
$\langle A\rangle^{-s}\cchi (H)P^{\perp}\langle A\rangle^{s}=\langle A\rangle^{-s}
[\cchi (H)P^{\perp},\langle A\rangle^{s}]+\cchi (H)P^{\perp}$. 
This is also true for $C\cchi (H)P^{\perp}C^{-1}$. 
By Theorem~\ref{th:mourre-reech-gene-proj}, \eqref{eq:r-lap-C} holds true. Since 
$E_\Ic(H)P^\perp=E_\Ic(H)$, $(H-z)^{-1}P^\perp=(H-z)^{-1}$ for $\re z\in\Ic '$. 
Therefore \eqref{eq:r-lap-C} yields \eqref{eq:lap}. \qed

\subsection{Energy estimates.} 
\label{ss:energy}

Here we extend a little bit G\'erard's method in \cite{ge}. We work in the general framework of 
Subsection~\ref{ss:special} and get the following improvements of 
Theorem~\ref{th:mourre-reech-gene-proj} and Theorem~\ref{th:mourre}. 
\begin{theorem}\label{th:mourre-reech-gene-proj-controle} 
Under the assumptions of Theorem~\ref{th:mourre-reech-gene-proj}, let 
$\sigma\in \{-1;1\}$ and 
choose a real $\mu$ such that $\sigma B'\geq 0$ with $B':=B+\mu$. Then 
\begin{align}
&\sup_{{\rm Re}z\in \Ic ', -\sigma {\rm Im}z>0}\|C(H-z)^{-1}P^\perp C\|
\label{eq:r-lap-C-controle} \\
&\hspace*{1cm}\leq \ 2\cdot \|CB'C^{-1}\|\cdot \|C^{-1}\cchi (H)P^{\perp}C\|\, +\, d^{-1}\cdot 
\|1-\cchi \|_{\rL^\infty}\cdot \|C\|^2\, , \nonumber
\end{align}
where $d$ is the distance between the support of 
$1-\cchi $ and $\Ic'$. 
\end{theorem} 
\begin{rem}\label{r:bound}
Note that, for $\sigma$ and $B$ as in Theorem~\ref{th:mourre-reech-gene-proj-controle}, 
one can always take $\mu =\sigma \|B\|$ to ensure $\sigma (B+\mu)\geq 0$. 
\end{rem}
\proof[Proof of Theorem~\ref{th:mourre-reech-gene-proj-controle}] By functional calculus, 
\begin{equation}\label{eq:energie-loin}  
\|C(H-z)^{-1}(1-\cchi )(H)P^\perp C\|\leq d^{-1}\|1-\cchi \|_{\rL^\infty}
\cdot \|C\|^2\, .
\end{equation}  
For $f\in\Hc$ and $z\in \C$ with $-\sigma{\rm Im}z>0$, let $u=(H-z)^{-1}\cchi (H)P^{\perp}Cf$. 
Notice that $E_\Ic (H)P^{\perp}u=u$. By \eqref{eq:mourre-reech-theta-chi} and a direct computation, 
\begin{align*}
\|Cu\|^2&\leq \langle u,[H,iB']u\rangle=2\im \langle B'u,(H-z)u\rangle
+2\sigma\im z\langle u \, ,\, \sigma B'u\rangle\\ 
&\leq 2\im \langle B'u,(H-z)u\rangle  \, ,
\end{align*} 
since $\sigma B'\geq 0$. Recall that 
$C\cchi (H)P^{\perp}C^{-1}$ is bounded. Thus $(H-z)u\in \Dc(C^{-1})$. In particular, since 
$CB'C^{-1}=CBC^{-1}+\mu$ is bounded, 
\begin{align}
2^{-1}\|Cu\|^2&\leq \im \langle CB'C^{-1} Cu,C^{-1}(H-z)u\rangle\nonumber
\leq \|CB'C^{-1}\|\cdot \|Cu\|\cdot \|C^{-1}(H-z)u\|,\nonumber
\end{align}
yielding 
$\|Cu\|\leq  2\|CB'C^{-1}\|\cdot\|C^{-1}\cchi (H)P^{\perp}C\|\cdot \|f\|$.
Together with \eqref{eq:energie-loin}, this implies \eqref{eq:r-lap-C-controle}.\qed

By combining the proof of Theorem~\ref{th:mourre} at the end of
Subsection~\ref{ss:special}  with
Theorem~\ref{th:mourre-reech-gene-proj-controle}, we derive:

\begin{corollary}\label{co:mourre-ameliore} 
Under the assumptions of Theorem~\ref{th:mourre}, take $s>1/2$,
$\sigma \in \{-1;1\}$, and  
$\cchi\in\Cc _c^\infty(\Ic )$ with $\cchi =1$ on $\Ic '$. Choose a real
$\mu$ such that  
$\sigma (\varphi (A)+\mu )\geq 0$. Then the l.h.s.\ of
\eqref{eq:lap} is bounded  
by the r.h.s.\  of \eqref{eq:r-lap-C-controle} for
$C=\sqrt{c_1}\langle A\rangle^{-s'}$ and $B'=\varphi (A)+\mu $ with  
$1/2<s'<1$, $2s'-1\leq \varepsilon _0$, and $s'\leq s$. 
\end{corollary} 
%

\subsection{Application.} 
\label{ss:application}

In practice, it is natural to try to derive a strict, weighted,
projected Mourre estimate  \eqref{eq:resc-mourre} from a similar
estimate containing some compact perturbation.  Precisely
\eqref{eq:resc-mourre} should follow from  
\begin{align}\label{eq:esti-mourre-proj-reech}  
& P^\perp E_\Ic (H)[H, i\varphi (A)]  E_\Ic(H) P^\perp \\
&\hspace*{1cm} \geq  E_\Ic(H)P^\perp \langle A\rangle ^{-(1+\varepsilon )/2}(
 c\, +\, K)\langle A\rangle ^{-(1+\varepsilon )/2}P^\perp E_\Ic(H), \nonumber
\end{align} 
for some compact operator $K$ and $c>0$. But to remove the influence
of $K$ using \eqref{eq:lim-fo}, we need to commute
$P^\perp E_\Ic(H)$ (or a regularized  version of it) through the
weight $\langle A\rangle ^{-(1+\varepsilon )/2}$. We are  able to do
this in the following situation.  
\begin{corollary}\label{co:mourre-reech-gene-proj} 
Let $\Ic$ be an open interval. Assume that, for all 
$\theta\in\Cc_c^\infty(\Ic ;\C)$, $P^\perp\theta (H)\in\Cc ^1(A)$. Let 
$\varepsilon_0\in ]0;1]$.
Assume further that, for all $\varepsilon \in ]0;\varepsilon_0]$, there exist $c>0$ and a 
compact operator $K$ such that, for all $R\geq 1$, there exists a real bounded 
borelian function $\varphi _R$ such that the weighted projected 
Mourre estimate 
\begin{align}\label{eq:esti-mou-proj-ree-R}  
& P^\perp E_\Ic (H)[H, i\varphi _R(A/R)]  E_\Ic(H) P^\perp \\
&\hspace*{1cm} \geq  P^\perp E_\Ic(H)\langle A/R\rangle ^{-(1+\varepsilon )/2}(
 c\, +\, K)\langle A/R\rangle ^{-(1+\varepsilon )/2}E_\Ic(H)P^\perp  \nonumber
\end{align} 
is satisfied. Then, for any $s>1/2$ and for 
any compact sub-interval $\Ic'$ of $\Ic$, the reduced LAP \eqref{eq:tal-reduit} for $H$ 
respectively to $(\Ic', s, A)$ holds true. 
\end{corollary} 
\proof By Remark~\ref{r:1/2<s<1}, we may assume that $s\in ]1/2;1[$ such that $\varepsilon :=2s-1
\in ]0;\varepsilon_0]$. By compactness of $\Ic'$, it is sufficient to show that, for any
$\lambda\in\Ic'$, \eqref{eq:tal-reduit} holds true with $\Ic'$ replaced by 
some open interval containing $\lambda$. It is enough to get \eqref{eq:tal-reduit} with $A$ replaced by $A/R$, for some $R\geq 1$. Let $\lambda\in\Ic'$. Since $K$ in 
\eqref{eq:esti-mou-proj-ree-R} is compact, we can use
\eqref{eq:lim-fo} to find $\cchi \in\Cc_c^\infty(\Ic ;\R)$ such that
$\cchi =1$ near $\lambda$ and $\|P^\perp\cchi (H)K\|\leq c/8$ (where $c$ appears 
in \eqref{eq:esti-mou-proj-ree-R}). Let $\Ic _1$ be an open
sub-interval of $\Ic'$ containing $\lambda$. From
\eqref{eq:esti-mou-proj-ree-R}, we get, for all $R\geq 1$, 
\begin{align}\label{eq:e-mou-pr-ree-theta}  
& \hspace{1.5cm}P^\perp E_{\Ic _1}(H)[H, i\varphi _R(A/R)] E_{\Ic
  _1}(H) P^\perp \ \geq \\ 
&  P^\perp E_{\Ic _1}(H)\langle A/R\rangle ^{-(1+\varepsilon )/2}\cdot \bigl(
 3c/4\, +\, (1-P^\perp\cchi (H))K(1-P^\perp\cchi (H))\bigr)\nonumber\\
&\hspace{2cm}\cdot \langle A/R\rangle ^{-(1+\varepsilon )/2}E_{\Ic _1}(H)P^\perp . 
\nonumber
\end{align} 
Since $1-P^\perp\cchi (H)=(1-\cchi )(H)+P\cchi (H)$, 
\begin{align*}
P^\perp E_{\Ic _1}(H)\langle A/R\rangle ^{-(1+\varepsilon
  )/2}(1-P^\perp\cchi (H))&=
 -P^\perp E_{\Ic _1}(H)\bigl[\langle A/R\rangle ^{-(1+\varepsilon )/2}\, ,\, 
P^\perp\cchi (H)\bigl]
\\
&=\, -P^\perp E_{\Ic _1}(H)\langle A/R\rangle ^{-(1+\varepsilon
  )/2}\cdot B_R,
\end{align*} 
where $\|B_R\|=O(1/R)$ by Propositions~\ref{p:regu} and~\ref{p:est3} (with $k=1$). 
Taking $R$ large enough (but fixed), we derive from
\eqref{eq:e-mou-pr-ree-theta}  the estimate 
\begin{eqnarray*}
P^\perp E_{\Ic _1}(H)[H, i\varphi _R(A/R)]E_{\Ic _1}(H) P^\perp 
\geq \ \frac{c}{2}\, P^\perp E_{\Ic _1}(H)\langle A/R\rangle ^{-1-\varepsilon }
E_{\Ic _1}(H)P^\perp . 
\end{eqnarray*} 
By Theorem~\ref{th:mourre-reech-gene-proj}, we obtain
\eqref{eq:r-lap-C} with $C=\langle A/R\rangle ^{-s}$ on some neighbourhood of $\lambda$,  
yielding \eqref{eq:tal-reduit} there with $A$ replaced by $A/R$.  \qed

\section{Perturbed Wigner-Von Neumann potentials.} 
\label{s:wigner}
\setcounter{equation}{0}

In this section, we apply our new theory to some special Schr\"odinger
operators (see Theorem~\ref{th:tal-Wigner}). As explained in
Section~\ref{s:intro}, we want to derive, on suitable intervals, a
usual Mourre estimate  (in Subsection~\ref{ss:wigner-mourre}) and a
weighted,  projected Mourre estimate (in
Subsection~\ref{ss:wigner-mourre-reech}) for the Schr\"odinger
operator $H_1$, see \eqref{e:H1}.

\subsection{Definitions and regularity.} 
\label{ss:def-regu}

Let $d\in\N^\ast$. We denote by $\langle \cdot ,\cdot \rangle$ and $\|\cdot\|$ the right 
linear scalar product and the norm in $\rL^2(\R^d)$, the space of squared integrable, complex 
functions on $\R^d$. 
Recall that $H_0$ is the self-adjoint realization of the Laplace operator 
$-\Delta$ in $\rL^2(\R^d)$ and that the Wigner-Von Neumann
potential $W:\R^d\longrightarrow \R$ is defined by $W(x)=q(\sin
k|x|)/|x|$, with $k>0$ and  $q\in\R^\ast$. Now we add to $W$ the
multiplication operator by the sum $V=V_{\rm sr}+V_{\rm lr}$  of real-valued
functions, $V_{\rm sr}$ has short range, and $V_{\rm lr}$ has long range.  
Precisely we require
\begin{assumption}\label{a:long-range} 
The functions $V_{\rm lr}$, $\langle x\rangle V_{\rm sr}$, and the distribution
$x\cdot \nabla V_{\rm lr}(x)$  belong to $\rL^\infty(\R^d)$. 
\end{assumption}
Under this assumption, on the Sobolev space $\rH^2(\R^d_x)$, the domain $\Dc (H_0)$ of $H_0$, 
\begin{eqnarray}\label{e:H1}
 H_1:=H_0+W+V = -\Delta + q |\cdot|^{-1}\sin (k|\cdot|)+V_{\rm sr}+V_{\rm lr}
\end{eqnarray} 
is self-adjoint. Let $P_1$ be the orthogonal
projection onto its pure point spectral subspace and
$P_1^\perp=1-P_1$. 
\\
Consider the strongly continuous one-parameter unitary group
$\{\Wc_t\}_{t\in \R}$ acting by:
\begin{equation}\label{groupe-c_0}
(\Wc_t f)(x)= e^{d t/2} f(e^tx), \mbox{ for all } f\in \rL ^2(\R^d_x).
\end{equation}
This is the $C_0$-group of dilation. A direct computation shows that 
\begin{equation}\label{domaine-conserve}
 \Wc_t\, \rH^2(\R^d_x) \subset \rH^2(\R^d_x), \mbox{ for
  all } t\in \R\, . 
\end{equation}
The generator of this group is the self-adjoint operator $A_1$, 
given by the closure of $(D_x\cdot x+ x\cdot D_x)/2$ on
$\Cc_c^\infty(\R^d)$ in $\rL^2(\R^d_x)$. For these reasons, the
operator $A_1$ is called \emph{the generator of dilation}. \\
The form $[W,iA_1]$ (defined on  $\Dc (A_1)\times\Dc
(A_1)$) extends to a bounded form associated to the multiplication
operator  by the function $W-W_1$, where 
\begin{equation}\label{eq:def-W_1}
 W_1(x)\ =\ qk\cos (k|x|)\ = \ (qk/2)\cdot (e^{ik|x|}\, +\, e^{-ik|x|})\, .
\end{equation}
In particular, $W\in \Cc^1(A_1)$ by
Proposition~\ref{p:caract-C1}. Furthermore, we prove 
\begin{proposition}\label{p:H_1-C1-A_1}  
We have $H_0\in \Cc^2(A_1)$. Moreover, under
Assumption~\ref{a:long-range}, the form $[V_{\rm   s}, iA_1]$, defined on 
$\Dc(A_1)\cap\Dc(H_0)$, extends to a bounded operator from $\rH^1(\R^d_x)$ to 
$\rH^{-1}(\R^d_x)$. In particular, $H_1\in \Cc^1(A_1)$. 
\end{proposition}  
\proof We use Section~\ref{s:usual}. As form on $\Dc(A_1)\cap
\Dc(H_0)$, $[H_0, iA_1]=2H_0$.  In particular, \eqref{e:C1b}
holds true with $A=A_1$ and $H=H_0$. By \eqref{domaine-conserve}   
and Theorem~\ref{th:abg2}, $H_0\in \Cc^1(A_1)$ and $[H_0, iA_1]_\circ
=2H_0$. For $z\not\in\R$, $R_0(z):=(H_0-z)^{-1}$ belongs to
$\Cc^1(A_1)$. Using \eqref{eq:comm-resolv} with $A=A_1$ and $H=H_0$,
we see that the form $[[R_0(z), iA_1]_\circ, iA_1]$ on $\Dc(A_1)\cap
\Dc(H_0)$ extends to bounded one.  Thus $R_0(z)\in \Cc^2(A_1)$ and
$H_0\in \Cc^2(A_1)$. \\ 
Since $\Dc(H_1)=\Dc(H_0)$ by Assumption~\ref{a:long-range}, $H_1\in
\Cc^1(A_1)$ follows from \eqref{domaine-conserve} and
Theorem~\ref{th:abg2} if \eqref{e:C1b} holds true with $A=A_1$ and  
$H=H_1$. We consider the form $[H_1, iA_1]$ on $\Dc(A_1)\cap
\Dc(H_0)$. It is the sum of $[H_0, iA_1]=2H_0$, of the bounded terms 
$[W,iA_1]=W-W_1$ and $[V_{\rm lr} ,iA_1]=x\cdot \nabla V_{\rm lr} $, and of 
\begin{align}
\langle f, [V_{\rm sr} ,iA_1]g\rangle &=\langle V_{\rm sr}f, iA_1g\rangle - \langle
A_1f, iV_{\rm sr}g\rangle \nonumber\\ 
&\label{e:commuVs} =\langle \langle x\rangle V_{\rm sr}f, \langle x\rangle
^{-1}(x\cdot \nabla 
_x+d/2)g\rangle 
+ 
\langle \langle x\rangle ^{-1}(x\cdot \nabla _x+d/2)f, \langle
x\rangle V_{\rm sr}g\rangle\, ,  
\end{align}
for $f,g\in \Dc(A_1)\cap \Dc(H_0)$. Since $\langle x\rangle V_{\rm sr}$ is bounded, 
\eqref{e:commuVs} extends to a bounded operator from $\rH^1(\R^d_x)$ to 
$\rH^{-1}(\R^d_x)$ and also from $\rH^2(\R^d_x)$ to $\rH^{-2}(\R^d_x)$. This gives 
\eqref{e:C1b} for $(H, A)=(H_1, A_1)$. Thanks to \eqref{domaine-conserve} and to 
Theorem~\ref{th:abg2}, $H_1\in \Cc^1(A_1)$. \qed

\subsection{Energy localization of oscillations.} 
\label{ss:oscillations}

\begin{figure}
\setlength{\unitlength}{1cm}
\begin{center}
\begin{picture}(8.5,6)(-1,-3)
%
\put(-1,0){\vector(1,0){8}} 
\put(0,-0.1){\line(0,1){0.2}}
\put(3,-0.1){\line(0,1){0.2}}
\put(6,-0.1){\line(0,1){0.2}}
\put(0,-0.5){0}
\put(2.7,-0.5){$k \hat x /2$}
\put(5.8,-0.5){$k \hat x$}
\put(0,0) {\arc(2,0){-110}}
\put(0,0) {\arc(2,0){100}}
\put(0,0) {\arc(2.5,0){-110}}
\put(0,0) {\arc(2.5,0){100}}
\put(6,0) {\arc(-2,0){110}}
\put(6,0) {\arc(-2,0){-100}}
\put(6,0) {\arc(-2.5,0){110}}
\put(6,0) {\arc(-2.5,0){-100}}
\put(-1,1){supp $\theta(|\cdot|^2)$ }
\put(4.6,1){supp $\theta(|\cdot- k\hat x|^2)$}
\end{picture}
\end{center}
\caption{supp $\theta\subset ]0, k^2/4[$.}\label{d:1}
\end{figure}
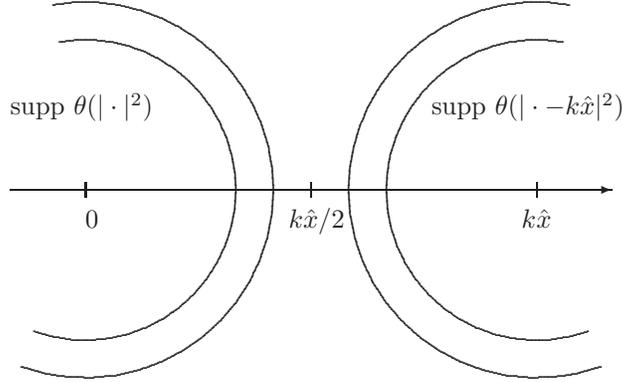

To prepare the derivation of Mourre estimates, we take advantage of
some ``smallness'' of  energy localizations of $W_1$ of the form
$\theta (H_0)W_1\theta (H_0)$, extending a result  by \cite{fh} in
dimension one. As seen in \cite{fh}, this term is not expected to be
small  if $\theta$ is localized near $k^2/4$. Using pseudodifferential
calculus, one have the same  
impression if $d\geq 2$ and if $\theta$ is supported in $]k^2/4 ;+\infty[$ (see 
Remark~\ref{r:above-threshold}). However, if $\theta$ lives in a small
enough compact  interval $\Ic \subset ]0;k^2/4[$, then the same
smallness as in \cite{fh} is valid as stated  in
Lemma~\ref{l:fh-local<seuil} below. The proof combines an idea in
\cite{fh} with pseudodifferential calculus (see
Subsection~\ref{ss:psi-do} for notation). In the sequel, we shall
write $\hat x$ for $x/|x|$.   
\begin{lemma}\label{l:fh-local<seuil}
Let $\lambda\in ]0;k^2/4[$. Recall that $g$ is given by
\eqref{eq:metric} and $W_1$ by  \eqref{eq:def-W_1}. Take $\cchi
_1\in\Cc ^\infty (\R^d)$ such that $\cchi _1=0$ near  $0$ and $\cchi
_1=1$ near infinity, and set $e_\pm(x)=\cchi _1(x)e^{\pm ik|x|}$.  For
$\theta\in \Cc _c^\infty(\R)$ with small enough support about
$\lambda$, there exist  symbols $b_0, b_{j, \sigma}\in \Sc (\langle
x\rangle ^{-1} \langle \xi\rangle ^{-1}, g)$, for $j\in\{1; 2\}$ and
$\sigma\in\{+, -\}$, such that  
\begin{equation}\label{eq:loca-oscill}
\theta (H_0)W_1\theta (H_0)\ =\ b_{1, +}^we_++b_{1, -}^we_-+\theta
(H_0)\bigl(e_+b_{2, +}^w +e_-b_{2, -}^w\bigr)+b_0^w\, .
\end{equation}
In particular, $\langle A_1\rangle^\varepsilon\theta (H_0)W_1\theta
(H_0)$ is compact on $\rL^2(\R^d_{x})$, for $\varepsilon\in [0;1[$.  
\end{lemma}
\begin{rem}\label{r:infini}
In dimension $d=1$, this result is proved in \cite{fh} and it also holds true if $\lambda >k^2/4$. 
Our proof below covers also this case. 
\end{rem}
\proof[Proof of Lemma~\ref{l:fh-local<seuil}] By pseudodifferential calculus, 
$\theta (H_0)(1-\cchi _1)W_1\theta (H_0)=b_0^w$ 
with $b_0\in\Sc (\langle x\rangle ^{-1}\langle \xi\rangle ^{-1}, g)$. By \eqref{eq:def-W_1} and 
the proof of Proposition~\ref{p:a^w-oscillant}, we can find $\cchi _3\in\Cc ^\infty (\R^d)$ 
such that $\cchi _3=0$ near $0$ and $\cchi _3=1$ near infinity, and $b_{j, \sigma}
\in\Sc (\langle x\rangle ^{-1}\langle \xi\rangle ^{-1}, g)$, for $j\in\{0;2\}$ 
and $\sigma\in\{+, -\}$, such that 
\begin{align}
&\hspace{-1cm}2(qk)^{-1}\theta (H_0)\cchi _1W_1\theta (H_0)\nonumber\\
&=\theta (H_0)\Bigl(\bigl(\theta (|\xi -k\hat{x}|^2)
\cchi _3(x)\bigr)^we_++\bigl(\theta (|\xi +k\hat{x}|^2)\cchi _3(x)\bigr)^we_-\Bigr)\nonumber\\
&\hspace{2.5cm}+\, \theta (H_0)\bigl(b_{0, +}^we_++b_{0, -}^we_-+e_+b_{2, +}^w+e_-b_{2, -}^w
\bigr)\nonumber\\
&=\bigl(\theta (|\xi |^2)\theta (|\xi -k\hat{x}|^2)
\cchi _3(x)\bigr)^we_++\bigl(\theta (|\xi |^2)\theta (|\xi +k\hat{x}|^2)\cchi _3(x)\bigr)^we_-
\nonumber\\
&\hspace{2.5cm}+\, b_{1, +}^we_++b_{1, -}^we_-+\theta (H_0)(e_+b_{2, +}^w+e_
-b_{2, -}^w)\, , \label{eq:dev-oscill}
\end{align} 
by composition. Now we choose the support of $\theta$ small enough about $\lambda$ such that $\theta (|\xi |^2)\theta (|\xi -k\hat{x}|^2)=0=\theta (|\xi |^2)
\theta (|\xi +k\hat{x}|^2)$, for all $x\neq 0$ and $\xi\in\R^d$. This
is possible since  $0\leq\lambda <k^2/4$, see Figure \ref{d:1}. Now \eqref{eq:dev-oscill} 
reduces to \eqref{eq:loca-oscill}.  
By Appendix~\ref{app:s:interpolation}, $\langle A_1\rangle^\varepsilon
\langle D_x\rangle^{-\varepsilon}\langle x\rangle^{-\varepsilon}$ extends to a bounded 
operator. For $b\in\Sc (\langle x\rangle ^{-1}\langle \xi\rangle ^{-1}, g)$, there exists 
$b\in\Sc (\langle x\rangle ^{\varepsilon -1}\langle \xi\rangle ^{\varepsilon -1}, g)$ such 
that $\langle x\rangle^{\varepsilon}\langle D_x\rangle^{\varepsilon}b^w=
b_1^w$ and $b_1^w$ is compact by \eqref{eq:caract-pseudo-compact}. Using 
\eqref{eq:loca-oscill}, this implies that $\langle A_1\rangle^\varepsilon\theta (H_0)W_1
\theta (H_0)$ is compact since we can write $\theta (H_0)e_+b_{2, +}^w=\theta (H_0)
\langle x\rangle ^{-1}e_+\langle x\rangle b_{2, +}^w$ 
with $\langle x\rangle b_{2, +}^w$ bounded and $\theta (H_0)\langle x\rangle ^{-1}=b^w$ with 
$b\in\Sc (\langle x\rangle ^{-1}\langle \xi\rangle ^{-1}, g)$. \qed

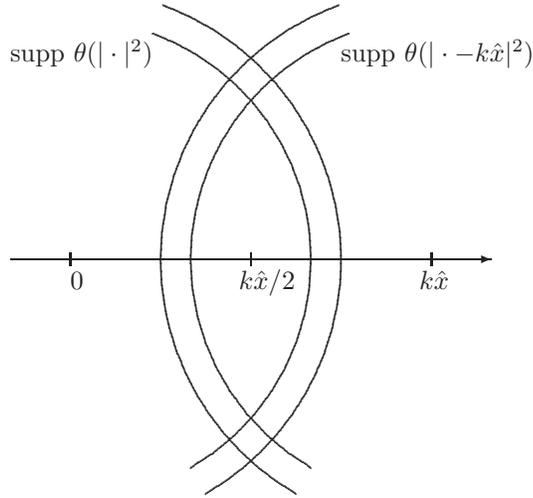
\begin{figure}
\setlength{\unitlength}{.8cm}
\begin{center}
\begin{picture}(8,8.5)(-1,-4)
\put(-1,0){\vector(1,0){8}} 
\put(0,-0.1){\line(0,1){0.2}} 
\put(3,-0.1){\line(0,1){0.2}}
\put(6,-0.1){\line(0,1){0.2}}
\put(0,-0.5){0}
\put(2.8,-0.5){$k \hat x /2$}
\put(5.8,-0.5){$k \hat x$}
\put(0,0) {\arc(4,0){-60}}
\put(0,0) {\arc(4,0){70}}
\put(0,0) {\arc(4.5,0){-60}}
\put(0,0) {\arc(4.5,0){70}}
\put(6,0) {\arc(-4,0){-70}}
\put(6,0) {\arc(-4,0){60}}
\put(6,0) {\arc(-4.5,0){-70}}
\put(6,0) {\arc(-4.5,0){60}}
\put(-1,3.3){supp $\theta(|\cdot|^2)$}
\put(4.5,3.3){supp $\theta(|\cdot- k\hat x|^2)$}
\end{picture}
\end{center}
\caption{supp $\theta\subset ]k^2/4, +\infty[$.}\label{d:2}
\end{figure}
\begin{rem}\label{r:above-threshold}
If $\lambda >k^2/4$ and $d>1$, the first two terms on the r.h.s.\  of
\eqref{eq:dev-oscill}  
do not vanish anymore, see Figure \ref{d:2}. In this case, our proofs
of the Mourre estimate  (cf., Proposition~\ref{p:mourre-hors-seuil})
and of the strict, weighted Mourre estimate  
(cf., Subsection~\ref{ss:wigner-mourre-reech}) do not work. \\
In dimension $d=1$, we note that the first two terms on the r.h.s.\
of \eqref{eq:dev-oscill} do vanish as soon as $\lambda\neq k^2/4$. 
See Figures \ref{d:1} and \ref{d:2} and recall that $\xi$ is co-linear
to $\hat x$. We recover a result in \cite{fh}.  
\end{rem}
%

\subsection{Usual Mourre estimate.} 
\label{ss:wigner-mourre}

Now we derive the Mourre estimate \eqref{eq:mourre} below $k^2/4$
under the following strengthening of Assumption~\ref{a:long-range}:  
\begin{assumption}\label{a:long-range-2}
The functions $V_{\rm lr}$, $\langle x\rangle V_{\rm sr}$, and the distribution
$x\cdot \nabla V_{\rm lr}(x)$  belong to $\rL^\infty(\R^d_x)$ and, as
operator of multiplication, compact from $\rH^2(\R^d_x)$ to $\rL^2(\R^d_x)$.
\end{assumption}
\begin{lemma}\label{l:diff-compacte}
Under Assumption~\ref{a:long-range-2}, $\varphi (H_1)-\varphi (H_0)$
is compact from $\rH^{2}(\R^d_x)$ to $\rL^{2}(\R^d_x)$, for $\varphi\in\Cc_c^\infty(\R)$.  
\end{lemma}
\proof 
Using \eqref{eq:int}, one has $\big(\varphi (H_1)-\varphi
(H_0)\big)\langle H_0\rangle=$ 
\begin{eqnarray*}
 \frac{i}{2\pi}\int_\C\partial_{\overline{z}}\,\varphi^\C(z)(z-H_1)^{-1}
 (W+V_{\rm sr}+V_{\rm lr})  \,   (z-H_0)^{-1}\langle H_0\rangle\, dz\wedge d\overline{z}.
\end{eqnarray*} 
For $z\notin \R$, the integrand is compact. Using \eqref{eq:majoA},
the integral converges in norm. Hence it is also compact.  \qed
\begin{proposition}\label{p:mourre-hors-seuil}
Under Assumption~\ref{a:long-range-2}, for any open interval $\Ic$ with 
$\overline{\Ic}\subset ]0;k^2/4[$, the Mourre estimate 
\eqref{eq:mourre} holds true for $(H, A)=(H_1, A_1)$. 
In particular, the point spectrum $\sigma _{pp}(H_1)$ of $H_1$ is
finite in $\Ic$. 
\end{proposition}
\proof It suffices to show \eqref{eq:mourre} on some compact
neighborhood of any $\lambda\in\Ic$. Take such a $\lambda\in\Ic$ and
let $\theta\in\Cc _c^\infty(\Ic;[0,1])$ such that $\theta =1$
near $\lambda$. Like in the proof of Proposition
\ref{p:H_1-C1-A_1}, as form on $\Dc (H_0)\cap \Dc(A_1) \times\Dc
(H_0)\cap \Dc(A_1)$, 
\[[H_0+V,iA_1]=2H_0- x\cdot \nabla V_{\rm lr}-\nabla \cdot x\langle x\rangle^{-1}
\langle x\rangle V_{\rm sr}-\langle x\rangle V_{\rm sr}\langle x\rangle^{-1}x\cdot
\nabla.\]
We recall that $[W,iA_1]_\circ=W-W_1$. Hence $[H_1,iA_1]_\circ$ is
bounded from $\Dc(H_0)=\rH^{2}(\R^d_x)$ to $\Dc(H_0)^*=\rH^{-2}(\R^d_x)$. Moreover, 
by Lemma~\ref{l:diff-compacte}, the bounded operator $\theta (H_1)[H_1,iA_1]_\circ\theta (H_1)$ 
is equal to $\theta (H_0)(2H_0-W_1)\theta(H_0)$, up to some compact operator. 
By Lemma~\ref{l:fh-local<seuil}, we can choose the support of $\theta$ such that $\theta 
(H_0)W_1\theta (H_0)$ is compact. Thus, there exist $c>0$ and compact operators  $K,K'$ 
such that 
\begin{eqnarray*}
 \theta (H_1)[H_1,A_1]_\circ\theta (H_1) \geq c\,\theta (H_0)^2\, +\, K'\,
 \geq c\,\theta (H_1)^2\, +\, K.  
\end{eqnarray*}
This yields the Mourre estimate \eqref{eq:mourre} near $\lambda$. \qed

As explained in Subsection~\ref{ss:reduced}, we need some information on
possible eigenvalues  embedded in the interval on which the LAP takes
place. Recall that $P_1$ denotes the orthogonal  projection onto the
pure point spectral subspace of $H_1$. 
\begin{proposition}\label{p:P_1-C^1}
Under Assumption~\ref{a:long-range-2}, take an open interval $\Ic$ with 
$\overline{\Ic}\subset ]0;k^2/4[$ such that, for all $\mu\in\Ic$,
${\rm Ker}(H_1-\mu)\subset 
\Dc (A_1)$. Then $E_{\Ic}(H_1)P_1\in\Cc ^1(A_1)$. 
\end{proposition}
\proof By Proposition~\ref{p:mourre-hors-seuil}, the point spectrum is finite in 
$\Ic$. Thus $E_{\Ic}(H_1)P_1\in\Cc ^1(A_1)$, by Proposition~\ref{p:rang-fini}.\qed

We now explain how to check the hypothesis ${\rm Ker}(H_1-\mu)\subset 
\Dc (A_1)$. The abstract Theorems given in \cite{c, fms1} do not apply
here because of the low regularity of $H_1$ w.r.t. $A_1$, see the inclusions
\eqref{eq:inclusions}, the implication \eqref{eq:global-impl-local}, and 
Proposition~\ref{p:pas-regu}. For $j\in\{1;\cdots ;d\}$, the multiplication operator by 
$x_j$ in $\rL^2(\R^d_x)$ is also denoted by $x_j$. As preparation, we show, using a 
Lithner-Agmon type equality, the following  
\begin{lemma}\label{l:Lithner-Agmon}
Let $n\in\N$. If $v\in\Cc^2(\R^d_x)\cap\rH ^2(\R^d_x)\cap\Dc (\langle x\rangle^{2n})$ then 
$\nabla _xv\in\Dc (\langle x\rangle^{n})$. 
\end{lemma}
\proof Define $\Phi (x)=n\ln \langle x\rangle$ for $x\in\R^d$ and let $R>1$. Using Green's 
formula, we can show that 
\begin{equation}\label{eq:lithner-agmon}
\int _{|x|\leq R}\bigl|\nabla (e^\Phi v)\bigr|^2dx=a(R)+\re\int _{|x|\leq R}e^{2\Phi}
\overline{v}\bigl(-\Delta v+|\nabla\Phi|^2v\bigr)dx\, ,
\end{equation}
where the term $a(R)$ contains surface integrals on $\{|x|=R\}$ and
tends to $0$ as  $R\to\infty$, thanks to $v\in\Dc (\langle
x\rangle^{2n})$ and $v\in\Cc ^2(\R^d)$.  Since $v\in\rH ^2(\R^d_x)$,
the last term in \eqref{eq:lithner-agmon} converges as  $R\to\infty$,
yielding $\nabla (e^\Phi v)\in\rL ^2(\R^d_x)$. Since $e^\Phi v\nabla\Phi 
\in\rL ^2(\R^d_x)$, $\langle x\rangle^{n}\nabla v=e^\Phi \nabla v\in\rL
^2(\R^d_x)$. \qed 
\begin{lemma}\label{l:regu-vect-propre}
Under Assumption~\ref{a:long-range-2}, let
$u\in\Cc^2(\R^d)\cap\rH^2(\R^d_x)$ and  $\lambda\in ]0;k^2/4[$ such
that $(H_1-\lambda )u=0$. Then $u\in\Dc (A_1)$.  Moreover, if $V_{\rm lr}
=0$, then $u=0$.  
\end{lemma}
\proof By Proposition~\ref{p:mourre-hors-seuil}, the usual Mourre
estimate holds true  near $\lambda$. Thus, one can apply Theorem 2.1
in \cite{fh}. Therefore  $u\in\Dc (\langle x\rangle^n)$, for all
$n\in\N$. By Lemma~\ref{l:Lithner-Agmon},  $D^\alpha u\in\Dc (\langle
x\rangle^{n})$, for all $n$ and all $\alpha\in\N^d$ with
$|\alpha|=1$. In particular, $x\cdot\nabla u\in\rL ^2(\R^d_x)$ and
$u\in\Dc (A_1)$.  If $V_{\rm lr} =0$, we can apply Theorem 14.7.2 in
\cite{h2} to $u$ yielding $u=0$. \qed 

\begin{rem}\label{r:regu-u}
If the potential $V=V_{\rm sr}+V_{\rm lr}$ belongs to $\Cc^m(\R^d)$ for some
integer $m>d/2$ then, by  elliptic regularity, any eigenvector $u$ of
$H_1$ belongs to $\Cc^2(\R^d)$. In particular,  by
Lemma~\ref{l:regu-vect-propre}, Proposition~\ref{p:P_1-C^1} applies to
any open interval $\Ic$ such that $\overline{\Ic}\subset ]0;k^2/4[$. 
\end{rem}
%

\subsection{Weighted Mourre estimate.} 
\label{ss:wigner-mourre-reech}

Here we establish for $H_1$ a projected, weighted Mourre estimate like 
\eqref{eq:esti-mou-proj-ree-R} in order to prove a limiting absorption principle (cf., 
Theorem~\ref{th:tal-Wigner}). 
To this end, we use the following assumption, which is stronger than 
Assumption~\ref{a:long-range-2}. 
\begin{assumption}\label{a:long-range-3}
For some $\rho_0 \in ]0,1]$, the functions $\langle x\rangle ^{\rho_0}V_{\rm lr}$, 
$\langle x\rangle ^{1+\rho_0}V_{\rm sr}$, and the distribution 
$\langle x\rangle ^{\rho_0}x\cdot \nabla V_{\rm lr}(x)$ belong to $\rL^\infty(\R^d_x)$. 
\end{assumption}

We start by strengthening Lemma \ref{l:diff-compacte}. 
\begin{lemma}\label{l:diff-compacte2}
Under Assumption~\ref{a:long-range-3}, for $\varepsilon\in [0;\rho
_0[$ and $\varphi\in\Cc_c^\infty(\R)$,  
\begin{equation}\label{eq:diff-funct-H_1-H_0-contre-A_1}
 (\varphi (H_1)-\varphi (H_0))\langle A_1\rangle ^{\varepsilon}\
 \mbox{is compact from }\rL^{2}(\R^d_x) \mbox{ to }\rH^{2}(\R^d_x).
\end{equation}
\end{lemma} 
\proof 
For $z\not\in\R$, $(z-H_0)^{-1}=r_z^w$ where $r_z$ satisfies
\eqref{eq:bound-r_z}  with $m=\langle \xi\rangle ^{2}$. By
composition, we can find, for all $\ell\in\N$,  $C_\ell>0$ and
$N_\ell\in\N$ such that, for all $z\not\in\R$,  
\[\|\langle x\rangle ^{-\rho _0}\# r_z\# (\langle x\rangle ^{\varepsilon}
\langle \xi\rangle ^{\varepsilon}) \|_{\ell , \Sc (\langle x\rangle
  ^{\varepsilon -\rho _0} \langle \xi\rangle ^{\varepsilon -2}, g)} 
\leq C_\ell \langle z\rangle^{N_\ell +1}|\im (z)|^{-N_\ell -1}\, .\]
Now thanks to Assumption~\ref{a:long-range-3},
\eqref{eq:caract-pseudo-compact}, \eqref{eq:dg1},  \eqref{eq:dg2}, and
\eqref{eq:dg3}, we infer that
\begin{align*}
\langle H_1\rangle \left(\varphi (H_1)-\varphi(H_0)\right) 
\langle x\rangle^{\varepsilon} 
\langle D_x\rangle^{\varepsilon}
&=\frac{i}{2\pi}\int_\C\partial_{\overline{z}}\,\varphi^\C(z) 
\langle H_1\rangle (z-H_1)^{-1} 
\\
& \hspace*{-1.5cm} (W+V_{\rm sr}+V_{\rm lr}) \langle x\rangle^{\rho _0} \langle
 x\rangle^{-\rho _0}   (z-H_0)^{-1}\langle x\rangle^{\varepsilon}
\langle D_x\rangle ^{\varepsilon}\, dz\wedge d\overline{z}
\end{align*}
is a compact operator, as norm convergent integral of compact
operators. To conclude, we recall $\langle x\rangle
^{-\varepsilon}\langle D_x\rangle ^{-\varepsilon} \langle A_1\rangle
^{\varepsilon}$ is bounded by Lemma~\ref{l:interpolation}. \qed

The main result on Schr\"odinger operators with
oscillating potential is 

\begin{theorem}\label{th:tal-Wigner}
Let $\lambda\in]0;k^2/4[$ and suppose that Assumption~\ref{a:long-range-3} is satisfied. 
Take a small enough, open interval $\Ic\subset ]0;k^2/4[$ about  $\lambda$ such that, for all
$\mu\in\Ic$, ${\rm Ker}(H_1-\mu)\subset\Dc (A_1)$. Then, for  any
$s>1/2$ and any interval $\Ic'\subset\overline{\Ic'}\subset\Ic$, the
reduced  LAP \eqref{eq:tal-reduit} for $H_1$ respectively to $(\Ic',
s, A_1)$ holds true.  
\end{theorem}
\begin{rem}\label{r:TAL-d=1->seuil}
Of course, a compactness argument shows that we can remove the
smallness condition on $\Ic$. We also get an estimate like in \eqref{eq:r-lap-C-controle} in 
Theorem~\ref{th:mourre-reech-gene-proj-controle}. \\
If $d=1$, the proof of Theorem~\ref{th:tal-Wigner} works also if $\lambda >k^2/4$, 
by Remark~\ref{r:infini}. \\
Under Assumption~\ref{a:long-range-3}, Theorem~\ref{th:tal-Wigner} ensures that 
$H_1$ has no singular spectrum above $\Ic$. In dimension $d=1$, for
$V_{\rm lr}=0$  but under a weaker assumption on $V_{\rm sr}$, this
result was already obtained in \cite{K, r} (see references
therein). Our long-range  result seems to be new, even in
dimension $1$.\\
If $\Ic$ contains an eigenvalue $\mu$ of $H_1$, the condition ${\rm Ker}(H_1-\mu)\subset
\Dc (A_1)$ is satisfied if $V=V_{\rm sr}+V_{\rm lr}$ is smooth enough (cf., Proposition~\ref{p:P_1-C^1} and 
Remark~\ref{r:regu-u}). 
\\
If one sets $q=0$, i.e.\ if one removes the potential $W$, $H_1$
  has no embedded eigenvalue (cf., \cite{fh}). The following proof of Theorem~\ref{th:tal-Wigner} 
works for each compact interval $\Ic '\subset (0, \infty)$ and gives the LAP \eqref{eq:lap} 
with $(H, A)=(H_1, A_1)$. This is a well-known result that can be obtained by the
 involved versions of the Mourre theory which are exposed in
\cite{abg,s}. But the technics of \cite{m,ge} do not apply, since the needed regularity condition 
$H_1\in \Cc^2(A_1)$ is not always satisfied under
Assumption~\ref{a:long-range-3}. 
\end{rem}
\proof[Proof of Theorem~\ref{th:tal-Wigner}:]
Let $\theta , \cchi , \tau\in \Cc _c^\infty  (]0;k^2/4[)$ such that
$\tau\cchi =\cchi$, $\cchi\theta=\theta$, and $\theta =1$ near 
$\Ic$. Later we shall adjust the size of the support of $\cchi$. 
By Proposition~\ref{p:H_1-C1-A_1} and \eqref{eq:global-impl-local}, $\cchi (H_1)\in
\Cc^1(A_1)$. Since $E_{\Ic}(H_1)P_1\in\Cc ^1(A_1)$ by Proposition~\ref{p:P_1-C^1}, 
$\cchi (H_1)P_1^\perp=\cchi (H_1)-\cchi (H_1)E_{\Ic}(H_1)P_1$ belongs to $\Cc^1(A_1)$. 

Let $s\in ]1/2;1[$. As in \cite{ge}, we define $\psi :\R\dans\R$ by 
\begin{equation}\label{def-phi}
\psi (t):=\int _{-\infty}^t\langle u\rangle ^{-2s}\, du . 
\end{equation}
Note that $\psi\in\Sc ^0$ and is in particular bounded. Let $R\geq 1$. 
As forms, using the fact that $H_1\tau (H_1)$  
is a bounded operator and belongs to $\Cc^1(A_1)$ and using \eqref{eq:egalite},  
\begin{align*}
 F&:=P_1^\perp\theta (H_1)[H_1,i\psi (A_1/R)]\theta
 (H_1)P_1^\perp=P_1^\perp\theta (H_1)[H_1\tau (H_1),i\psi
 (A_1/R)]\theta (H_1)P_1^\perp\\ 
&=\frac{i}{2\pi}\int_\C\partial_{\overline{z}}\psi^\C(z)P_1^\perp\theta (H_1)
(z-A_1/R)^{-1}[H_1\tau (H_1),iA_1/R]_\circ\\
&\hspace{7cm}(z-A_1/R)^{-1}\theta (H_1)P_1^\perp dz\wedge 
d\overline{z}\, .
\end{align*}
Next to $P_1^\perp\theta (H_1)$ we let appear $\cchi
(H_1)P_1^\perp$ and commute it with $(z-A_1/R)^{-1}$.
Since $\cchi (H_1)P_1^\perp\in\Cc^1(A_1)$, 
%
%
we obtain, using \eqref{eq:dg1}, \eqref{eq:dg2}, and \eqref{eq:majoA}, for some uniformly 
bounded operator $B_1$ w.r.t. $R\geq 1$, 
\begin{align}
F=&\frac{i}{2\pi}\int_\C\partial_{\overline{z}}\psi^\C(z)P_1^\perp\theta (H_1)
(z-A_1/R)^{-1}P_1^\perp\cchi (H_1)[H_1\tau (H_1),iA_1/R]_\circ\nonumber\\
&\hspace{5cm}\cchi (H_1)P_1^\perp (z-A_1/R)^{-1}\theta (H_1)P_1^\perp dz\wedge 
d\overline{z}\label{eq:F}\\
&+\ P_1^\perp\theta (H_1)\langle A_1/R\rangle ^{-s}R^{-2}B_1\langle A_1/R\rangle ^{-s}
\theta (H_1)P_1^\perp\, .\nonumber 
\end{align}
Let $\varepsilon:=\rho _0/2$. Using
\eqref{eq:diff-funct-H_1-H_0-contre-A_1}, notice that 
\begin{align*}
G&:=P_1^\perp\cchi (H_1)[H_1\tau (H_1),iA_1/R]_\circ\cchi
(H_1)P_1^\perp=P_1^\perp\cchi (H_1)[H_1,iA_1/R]_\circ\cchi
(H_1)P_1^\perp\\ 
&=P_1^\perp\cchi (H_1)[H_1,iA_1/R]_\circ\cchi (H_0) P_1^\perp +
P_1^\perp\cchi (H_1)K_1R^{-1}B_2\langle A_1/R\rangle
^{-\varepsilon}P_1^\perp,  
\end{align*}
where the operator $K_1:=\tau (H_1)[H_1,iA_1]_\circ(\cchi (H_1)-\cchi
(H_0))\langle 
A_1\rangle ^{\varepsilon}$  is compact and $B_2:=\langle A_1/R\rangle
^{\varepsilon}\langle A_1\rangle ^{-\varepsilon}$ is  uniformly
bounded. Similarly, there is $K_2$ compact so that 
\begin{align}
G=&P_1^\perp\cchi (H_0)[H_1,iA_1/R]_\circ\cchi (H_0)P_1^\perp
+P_1^\perp\cchi (H_1)K_1R^{-1} B_2\langle A_1/R\rangle
^{-\varepsilon}P_1^\perp\nonumber
\\ &+\, P_1^\perp \langle
A_1/R\rangle ^{-\varepsilon}B_2K_2R^{-1}\cchi (H_0)P_1^\perp.\label{eq:G} 
\end{align}

We focus on the potential contribution in $G$. 
Choosing $\tau$ appropriately and using Assumption~\ref{a:long-range-3}, we claim that there exist a 
compact operator $K_3$ and an uniformly bounded operator $B_3$ such that  
\begin{align}\label{e:potentialpart}
\cchi (H_0)[W+V,iA_1/R]_\circ\cchi (H_0)&=R^{-1}\cchi (H_0)K_3B_3\langle
A_1/R\rangle ^{-\varepsilon} \cchi (H_0).
\end{align}
By Lemma~\ref{l:fh-local<seuil}, we take the
support of $\tau$ small enough to ensure the compactness of $\tau (H_0)  W_1\tau (H_0)
\langle A_1\rangle ^{\varepsilon}$. By writing   
\[(W+ x\cdot \nabla V_{\rm lr})\tau (H_0)\langle A_1\rangle
^{\varepsilon} = (W+ x\cdot \nabla V_{\rm lr})\langle x\rangle^{\rho _0} \cdot  
b^w\cdot \langle x\rangle ^{-\varepsilon}\langle D_x\rangle
^{-\varepsilon}\langle A_1\rangle ^{\varepsilon},\] 
with $b\in\Sc (\langle x\rangle ^{-\varepsilon }\langle \xi\rangle
^{-1}, g)$, $\tau (H_0)[W+ V_{\rm lr},iA_1]_\circ\tau (H_0)\langle A_1\rangle
^{\varepsilon}$ is compact by \eqref{eq:caract-pseudo-compact}  and
Lemma~\ref{l:interpolation}. Similarly, we prove the compactness of 
$\tau(H_0)[V_{\rm sr},iA_1]_\circ\tau (H_0)\langle A_1\rangle ^{\varepsilon}$, 
making use of the fact that, by \eqref{e:commuVs}, $\langle D_s\rangle ^{-1}
[V_{\rm sr},iA_1]_\circ\langle x\rangle^{\rho _0}\langle D_s\rangle ^{-1}$ extends 
to a bounded operator. This yields \eqref{e:potentialpart}.\\
Taking advantage of  $[H_0,iA_1]_\circ=2H_0$, of \eqref{e:potentialpart}, and of
\eqref{eq:G}, we rewrite \eqref{eq:F}:
\begin{align}
F&=\frac{i}{2\pi}\int_\C\partial_{\overline{z}}\psi^\C(z)P_1^\perp\theta (H_1)
(z-A_1/R)^{-1}P_1^\perp2R^{-1}H_0\cchi ^2(H_0)\nonumber\\
&\hspace{6cm}P_1^\perp (z-A_1/R)^{-1}\theta (H_1)P_1^\perp dz\wedge 
d\overline{z}\label{eq:F1}\\
&\hspace*{1cm}+\ P_1^\perp\theta (H_1)\langle A_1/R\rangle
^{-s}(R^{-2}B_1+R^{-1}K_4) 
\langle A_1/R\rangle ^{-s}
\theta (H_1)P_1^\perp\, , \nonumber
\end{align}
with compact $K_4$ such that, for some $c_1>0$, 
\begin{equation}\label{eq:bound-K}
\|K_4\|\, \leq \, c_1(\|P_1^\perp\cchi (H_1)K_1\|+\|K_2\cchi (H_0)\|+\|\cchi (H_0)K_3\|)\, . 
\end{equation}
Next we commute $(z-A_1/R)^{-1}$ with $P_1^\perp2R^{-1}H_0\cchi
^2(H_0)P_1^\perp$. Recalling \eqref{eq:int} with $k=1$ and \eqref{def-phi}, there are 
$B_4$ and $B_5$, uniformly bounded, such that
\begin{align*}
F&=P_1^\perp\theta (H_1)\psi '(A_1/R)P_1^\perp2R^{-1}H_0\cchi ^2(H_0)P_1^\perp
\theta (H_1)P_1^\perp\\
&\hspace*{1cm} +\ P_1^\perp\theta (H_1)\langle A_1/R\rangle
^{-s}(R^{-2}B_4+R^{-1}K_4) 
\langle A_1/R\rangle ^{-s}
\theta (H_1)P_1^\perp\, , 
\\
&=P_1^\perp\theta (H_1)\langle A_1/R\rangle ^{-s}2R^{-1}H_0\cchi ^2(H_0)
\langle A_1/R\rangle ^{-s}\theta (H_1)P_1^\perp\\
&\hspace*{1cm} +\ P_1^\perp\theta (H_1)\langle A_1/R\rangle ^{-s}(R^{-2}B_5+R^{-1}K_4)
\langle A_1/R\rangle ^{-s}
\theta (H_1)P_1^\perp\, , 
\\
&\geq 2R^{-1}c_2P_1^\perp\theta (H_1)\langle A_1/R\rangle ^{-s}\cchi ^2(H_0)
\langle A_1/R\rangle ^{-s}\theta (H_1)P_1^\perp\\
&\hspace*{1cm} +\ P_1^\perp\theta (H_1)\langle A_1/R\rangle
^{-s}(R^{-2}B_5+R^{-1}K_4) 
\langle A_1/R\rangle ^{-s}
\theta (H_1)P_1^\perp\, , 
\end{align*}
where $c_2>0$ is the infimum of
$\Ic$. Finally, since $K_5:=\cchi ^2(H_0)- \cchi ^2(H_1)$ is compact
by  \eqref{eq:diff-funct-H_1-H_0-contre-A_1}, we find an uniformly
  bounded $B_6$, such that
\begin{align*}
F&\geq 2R^{-1}c_2P_1^\perp\theta (H_1)\langle A_1/R\rangle ^{-s}\cchi ^2(H_1)
\langle A_1/R\rangle ^{-s}\theta (H_1)P_1^\perp\\
&\hspace*{1cm}+\, P_1^\perp\theta (H_1)\langle A_1/R\rangle
^{-s}(R^{-2}B_5+R^{-1}K_4+ 
R^{-1}K_5)\langle A_1/R\rangle ^{-s}\theta (H_1)P_1^\perp\\
&\geq 2R^{-1}c_2 P_1^\perp\theta (H_1)\langle A_1/R\rangle
^{-2s}\theta (H_1)P_1^\perp  
+\ P_1^\perp\theta (H_1)\langle A_1/R\rangle ^{-s}\cdot \\
&\hspace*{1cm}(R^{-2}B_6+R^{-1}K_4+R^{-1}K_5\cchi (H_1)P_1^\perp)
\langle A_1/R\rangle ^{-s}\theta (H_1)P_1^\perp \, . 
\end{align*}
To conclude, using \eqref{eq:bound-K},  we decrease the support of $\cchi$ to
ensure that $\|K_4\|+ \|K_5\cchi (H_1)P_1^\perp \|< c_2$. Subsequently, we
choose  $R>1$ large enough to guarantee $F\geq R^{-1}c_2P_1^\perp\theta
(H_1) \langle A_1/R\rangle ^{-2s}\theta (H_1)P_1^\perp$.  Letting act
the projector $E_\Ic (H_1)$ on both sides of this inequality and recalling  
the definition of $F$, we get the projected, weighted Mourre estimate 
\eqref{eq:mourre-reech-theta-chi} with $H=H_1$, $P=P_1$, $B=\psi
(A_1/R)$, and $C= \sqrt{c_2/R}\langle A_1/R\rangle ^{-s}$. By
Theorem~\ref{th:mourre-reech-gene-proj},  
we obtain the result. \qed
%

%

\section{Usual Mourre theory.} \label{s:usual}
\setcounter{equation}{0}

%
In this section, we explain why the usual Mourre theory with conjugate
operator $A_1$  cannot be applied to $H_1$, the considered Schr\"odinger operator
with oscillating potential.  We have proved that $H_1\in \Cc^1(A_1)$
and established a Mourre estimate for $H_1$
w.r.t.\ $A_1$, see  Propositions~\ref{p:H_1-C1-A_1} and
\ref{p:mourre-hors-seuil}. However, in order to apply the standard
Mourre theory, one has to prove that $H_1$ is in a better class of
regularity w.r.t.\ $A_1$. In this section, we prove that this is not
the case. If one replaces $A_1$ by some natural variants, we explain in 
Remark~\ref{rq:autre-op-conj} below that the required regularity is not
available. On the other hand, a consequence of Theorem
\ref{th:tal-Wigner} is that under the Assumption \ref{a:long-range-3},
the operator $H_1$ has no singular continuous spectrum. By abstract
means, see \cite[Proposition 7.2.14]{abg}, there exists a conjugate
operator $\tilde A$, such that $H_1\in \Cc^\infty(\tilde A)$ and such
that a strict Mourre estimate holds true for $H_1$, w.r.t.\ $\tilde
A$, on every interval that contains neither an eigenvalue nor $\{0,
k^2/4\}$. It seems very difficult to find explicitly $\tilde A$.

We first continue the description of different classes of regularity
appearing in the Mourre theory that we began in
Subsection~\ref{ss:regularity}. We refer again to  \cite{abg,ggm1,gg}
for more details. Recall that a self-adjoint operator $H$ belongs to
the  class $\Cc^1(A)$ if, for some (hence for all) $z\notin
\sigma(H)$, the bounded operator $(H-z)^{-1}$ belongs to
$\Cc^1(A)$. Lemma~6.2.9 and Theorem~6.2.10 in \cite{abg} gives the
following  characterization of this regularity:  
\begin{theorem}(\cite{abg})\label{th:abg} 
Let $A$ and $H$ be two self-adjoint operators in the Hilbert space
$\Hr$. For $z\notin \sigma(H)$, set $R(z):=(H-z)^{-1}$. The following
points are equivalent:  
\begin{enumerate} 
\item $H\in\Cc^1(A)$. 
\item For one (then for all) $z\notin \sigma(H)$, there is a finite
$c$ such that 
\begin{align}\label{e:C1a} 
|\langle A f, R(z) f\rangle - \langle R(\overline{z}) f, Af\rangle| \leq c 
\|f\|^2, \mbox{ for all $f\in\Dc(A)$}. 
\end{align} 
\item 
\begin{enumerate} 
\item [a.]There is a finite $c$ such that for all $f\in \Dc(A)\cap\Dc(H)$: 
\begin{equation}\label{e:C1b} 
|\langle Af, H f\rangle- \langle H f, Af\rangle|\leq \,
 c\big(\|H f\|^2+\|f\|^2\big). 
\end{equation} 
\item [b.] The set
$\{f\in\Dc(A);\,  R(z)f\in\Dc(A)\, \mbox{\rm and}\, R(\overline{z})f\in\Dc(A)
\}$ is a core for $A$, for some (then for all) $z\notin \sigma(H)$. 
\end{enumerate} 
\end{enumerate} 
\end{theorem} 
Note that the condition (3.b) could be uneasy to check, see
\cite{gg}. We mention \cite{GoleniaMoroianu}[Lemma A.2]
to overcome this subtlety. Note that \eqref{e:C1a} yields that the
commutator $[A, R(z)]$ extends to a bounded operator, in the form
sense. We shall denote the extension by $[A, R(z)]_\circ$. In the same
way, from \eqref{e:C1b}, the commutator $[H, A]$ extends to a unique
element of $\Bc\big(\Dc(H), \Dc(H)^*\big)$ denoted by $[H,
  A]_\circ$. Moreover, if $H\in \Cc^1(A)$ and $z\notin \sigma(H)$, 
\begin{eqnarray}\label{eq:comm-resolv}
\big[A, (H-z)^{-1}\big]_\circ =\quad  \underbrace{(H-z)^{-1}}_{\Hr
  \leftarrow \Dc(H)^*}\quad  \underbrace{[H, A]_\circ}_{\Dc(H)^*\leftarrow
  \Dc(H)} \quad \underbrace{(H-z)^{-1}}_{\Dc(H)\leftarrow \Hr}.
\end{eqnarray} 
Here we use the Riesz lemma to identify $\Hr$ with its anti-dual
$\Hr^*$. It turns out that an easier characterization is available if 
the domain of $H$ is conserved under the action of the unitary group 
generated by $A$. 
\begin{theorem}(\cite[p.\ 258]{abg})\label{th:abg2} 
Let $A$ and $H$ be two self-adjoint operators in the Hilbert space
$\Hr$ such that $e^{itA}\Dc (H)\subset\Dc (H)$, for all $t\in\R$. 
Then $H\in\Cc^1(A)$ if and only if \eqref{e:C1b} holds true. 
\end{theorem}
\begin{rem}\label{r:C1dom}
Some arguments used in the proof of Proposition~\ref{p:H_1-C1-A_1} may be performed 
in an abstract way. Take a Hilbert space $\Gr$, such that
$\Gr\hookrightarrow \Hr$ with a continuous, dense embedding and such that 
the $C_0$-group $\{e^{itA}\}_{t\in \R}$ stabilizes $\Gr$ (hence also $\Gr^*$ by duality). Let 
$T\in\Bc(\Gr, \Gr^*)$. We say that $T\in \Cc^1(A; \Gr, \Gr^*)$ if the strong limit of
$t\mapsto t^{-1}(e^{itA}Te^{-itA}-T)$ exists in $\Bc(\Gr, \Gr^*)$, as $t$ goes
to $0$. The limit is denoted by $[T, iA]_\circ$. Assuming the invariance of $\Dc (H)$ under 
the $C_0$-group $\{e^{itA}\}_{t\in \R}$ and taking a $\Gr$ with continuous, dense embeddings 
$\Dc (H)\hookrightarrow\Gr\hookrightarrow\Hr$, then $\{e^{itA}\}_{t\in \R}$ stabilizes $\Gr$ 
by interpolation. If $T\in \Cc^1(A; \Gr, \Gr^*)$ then $T\in\Cc^1(A; \Dc(H), \Dc(H)^*)$. 
If $H\in\Cc^1(A; \Dc(H), \Dc(H)^*)$, $[H, iA]_\circ$ coincide with the previous definition.
One can reformulate Theorem \ref{th:abg2} as follows: $H\in \Cc^1(A; \Dc(H), \Dc(H)^*)$ if 
and only if $H\in\Cc^1(A)$. 
\end{rem}   

We need to introduce others classes inside $\Cc^1(A)$. Let $T\in \Bc (\Hr)$. We say that 
$T\in \Cc^{1, u}(A)$ if the map 
$\R\ni t\mapsto e^{itA}Te^{-itA}\in \Bc (\Hr )$ has the usual $\Cc^1$ regularity.
We say that $T\in\Cc^{1,1}(A)$ if 
\begin{equation}\label{eq:def-c1,1}
\int_0^1 \big\|[[T, e^{itA}], e^{itA}]\big\|\, t^{-2}\, dt\ <\ \infty.
\end{equation}
We say that $T\in\Cc^{1+0}(A)$ if $T\in \Cc^1(A)$ and 
\begin{equation}\label{eq:def-c1+0}
\int_0^1 \big\|e^{itA}[T,A]e^{-itA}\big\|\, t^{-1}\, dt\ <\ \infty.
\end{equation}
Thanks to \cite[p.\ 205]{abg}, it turns out that 
\begin{equation}\label{eq:inclusions}
\Cc^2(A)\subset \Cc^{1+0}(A)\subset \Cc^{1,1}(A)\subset \Cc^{1, u}(A)\subset \Cc^1(A).
\end{equation} 
Given a self-adjoint operator $H$ and an open interval $\Ic$ of $\R$, we consider 
the corresponding local classes defined by: $H\in\Cc^{[\cdot]}_\Ic(A)$ if, for 
all $\varphi\in\Cc^\infty_c(\Ic)$, $\varphi (H)\in\Cc^{[\cdot]}(A)$. 
We say that $H\in\Cc^{[\cdot]}(A)$ if, for some $z\not\in\sigma (H)$, 
$R(z)\in\Cc^{[\cdot]}(A)$. Proposition~\ref{p:global-impl-local} also 
works for the new classes: for all open interval $\Ic$ of $\R$ and all $\varphi\in
\Cc^\infty _c(\Ic)$, 
\begin{equation}\label{eq:global-impl-local}
H\in\Cc^{[\cdot]}(A)\hspace{.4cm}\impl\hspace{.4cm}\varphi (H)\in\Cc^{[\cdot]}_\Ic(A)\, .
\end{equation}
In \cite{abg}, the LAP is obtained for $H\in\Cc^{1,1}(A)$ (see p.\ 308 and p.\ 317) and this
class is shown  to be optimal among the global classes (see the end of
Section 7.B). In \cite{s}, for  $H\in\Cc^{1+0}_\Ic(A)$, the LAP is
obtained on compact sub-interval of $\Ic$. It is  expected that the
class $\Cc^{1,1}_\Ic(A)$ is sufficient. Section 7.B in \cite{abg}
again shows that one cannot use in general a bigger local class to get
the LAP.   \\
Now we explore the regularity properties of $H_1$ under
Assumption~\ref{a:long-range-2}.  From Proposition~\ref{p:H_1-C1-A_1},
we know that $H_1\in\Cc^{1}(A_1)$. If $H_1$ would belong to
$\Cc^{1,1}(A_1)$ then, by \eqref{eq:inclusions} and
\eqref{eq:global-impl-local},  $H_1$ would belong to $\Cc^{1, 
  u}_{\Ic}(A_1)$ for any open interval $\Ic\subset]0;+\infty[$.  If
$H_1$ would belong to $\Cc^{1+0}_\Ic (A_1)$ or even to $\Cc^{1,1}_\Ic
(A_1)$, for  some open interval $\Ic\subset]0;+\infty[$, then $H_1$
would belong to $\Cc^{1,   u}_{\Ic}(A_1)$  by
\eqref{eq:inclusions}. In both cases, this would contradict:  
\begin{proposition}\label{p:pas-regu}
Under Assumption~\ref{a:long-range-2}, for any open sub-interval $\Ic$
of $]0;+\infty[$,  $H_1\not\in\Cc^{1, u}_{\Ic}(A_1)$. 
\end{proposition} 
\proof Take such an interval $\Ic$ and $\varphi\in\Cc^\infty _c(\Ic)$. 
By Proposition~\ref{p:H_1-C1-A_1}, \eqref{eq:inclusions}, and \eqref{eq:global-impl-local}, 
$\varphi (H_0)\in\Cc^{1, u}(A_1)$. Assume that $\varphi (H_1)\in\Cc^{1, u}(A_1)$. 
Then $K:=\varphi (H_1)-\varphi (H_0)\in\Cc^{1, u}(A_1)$ and $K$ is a
compact operator on  $\rL^2(\R^d_x)$, thanks to
Lemma~\ref{l:diff-compacte}. Thus  
\[[K,iA_1]_\circ =\lim _{t\to 0}t^{-1}\bigl(e^{-itA_1}Ke^{itA_1}-K\bigr)\]
in $\Bc (\rL^2(\R^d_x))$ and $[K,iA_1]_\circ$ is also compact. So is
$B[K,iA_1]_\circ B'$, for any  
$B, B'\in\Bc (\rL^2(\R^d_x))$. This contradicts
Lemma~\ref{l:non-compact} below. \qed 

\begin{lemma}\label{l:non-compact}
Assume Assumption~\ref{a:long-range-2}. For any open interval
$\Ic\subset ]0;+\infty[$, there exist  a function
$\varphi\in\Cc^\infty _c(\Ic)$ and bounded operators $B, B'$ on
$\rL^2(\R^d_x)$  such that  $B[\varphi (H_1)-\varphi (H_0),iA_1]_\circ B'$ is
not compact on $\rL^2(\R^d_x)$.  
\end{lemma}
We refer to Appendix \ref{app:s:dim1} for a proof of this Lemma for
$d=1$, which does not rely on pseudodifferential calculus.

\proof[Proof of Lemma \ref{l:non-compact}] 
In the sequel, for $C,D\in\Bc(\rL^2(\R^d_x))$, we write
$C\simeq D$ if $C-D$ is compact on $\rL^2(\R^d_x)$.
By Proposition~\ref{p:H_1-C1-A_1}, $H_1, H_0\in \Cc^1(A)$. Then $B_1:=
[\varphi (H_1)-\varphi (H_0),iA_1]_\circ$ is bounded. Furthermore, thanks to
\eqref{eq:egalite}, \eqref{eq:reste22}, and by the resolvent formula,
with a norm convergent integral, 
\begin{align}
 B_1 =\frac{i}{2\pi}\int_\C\partial _{\overline{z}}\,\varphi^\C
(z)\bigl[(z-H_1)^{-1}(W+V)(z-H_0)^{-1} 
, iA_1\bigr]_\circ\, dz\wedge d\overline{z}. \label{eq:B_1=}
\end{align}
We recall that given a continuous function $F:\R^d\rightarrow \C$,
that tends to $0$ at infinity, the multiplication by $F$ is compact
from $\rH^s(\R^d_x)$ to $L^2(\R^d_x)$, for all $s>0$. Using 
again Proposition~\ref{p:H_1-C1-A_1} and expanding the commutator, 
using the computation of $[W,iA_1]_\circ$ (see just before
\eqref{eq:def-W_1}) and  the resolvent formula again, it yields:
\begin{align}
 B_1 &\simeq\ \frac{i}{2\pi}\int_\C\partial _{\overline{z}}\varphi^\C
 (z)(z-H_1)^{-1} 
[W+V,iA_1]_\circ(z-H_0)^{-1}\, dz\wedge d\overline{z},
\nonumber
\\
&\simeq \frac{-i}{2\pi}\int_\C\partial _{\overline{z}}\varphi^\C (z)(z-H_1)^{-1}W_1
(z-H_0)^{-1}\, dz\wedge d\overline{z},
\nonumber
\\
\label{eq:B1-equiv}
&\simeq\ \frac{-i}{2\pi}\int_\C\partial _{\overline{z}}\varphi^\C (z)(z-H_0)^{-1}W_1
(z-H_0)^{-1}\, dz\wedge d\overline{z},
\\
&\label{eq:B1-equiv-loin-0}
\simeq\ \frac{-i}{2\pi}\int_\C\partial _{\overline{z}}\varphi^\C
(z)(z-H_0)^{-1}\cchi _1W_1 
(z-H_0)^{-1}\, dz\wedge d\overline{z},
\end{align}
with $\cchi _1\in\Cc ^\infty(\R^d)$, $\cchi _1=0$ near $0$, and $\cchi
_1=1$ near infinity. 

At this point, we use pseudodifferential techniques and, in particular, 
Appendix~\ref{app:s:oscillation}. 
For $x\in\R^d$, let $e_\pm (x)=\cchi _1(x)e^{\pm ik|x|}$. 
By \eqref{eq:def-W_1}, $(\cchi
_1W_1)(x)=kq2^{-1}(e_+(x)+e_-(x))$.  
Now we apply Proposition~\ref{p:a^w-oscillant} to $a(x, \xi)=|\xi
|^2\in\Sc (\langle \xi\rangle ^2,  
g)$. By its proof, $a_\pm$ can be chosen real and $a_\pm^w$ 
is self-adjoint. Using the resolvents of $a_\pm^w$ and $a^w$,  
\begin{align}
e_\pm(z-H_0)^{-1}&=e_\pm(z-a^w)^{-1}\ =\ (z-a_\pm^w)^{-1}e_\pm\nonumber\\
&\hspace*{1cm}+(z-a_\pm^w)^{-1}(e_\pm b_\pm^w+c_\pm^we_\pm)(z-a^w)^{-1}\,
, \label{eq:resolv-oscillant} 
\end{align}
for all $z\not\in\R$. 
We obtain from \eqref{eq:B1-equiv-loin-0}, 
\eqref{eq:resolv-oscillant}, and Proposition~\ref{p:a^w-oscillant}: 
\begin{equation*}
B_1\ \simeq\ -\sum _{\sigma =\pm}\frac{ikq}{4\pi}\int_\C
\partial _{\overline{z}}\varphi^\C (z)(z-H_0)^{-1}(z-a_\sigma ^w)^{-1}e_\sigma\, 
dz\wedge d\overline{z}\, .
\end{equation*}
According to \cite{bo} (see Appendix~\ref{app:s:calcul-fonct-pseudo}), $(z-H_0)^{-1}=p_z^w$ 
and $(z-a_\sigma ^w)^{-1}=p_{\sigma , z}^w$ where the symbols $p_z, p_{\sigma , z}$ belong 
to $\Sc (\langle \xi\rangle^{-2}, g)$ and satisfy \eqref{eq:bound-r_z} with $m=
\langle \xi\rangle^{2}$. Using the continuity of the map $\Sc (\langle \xi\rangle^{-2}, g)^2
\ni (r, t)\donne r\# t -rt\in\Sc (h\langle \xi\rangle^{-2}, g)$, we can find, for all 
$\ell\in\N$, $C_\ell>0$ and $\N_\ell\in\N$ such that
\begin{equation}\label{eq:esti-z}
\|p_z\# p_{\sigma , z}-p_zp_{\sigma , z}\|_{\ell , \Sc (h\langle \xi\rangle^{-4}, g)}\leq 
C_\ell \langle z\rangle^{N_\ell +1}
|\im (z)|^{-N_\ell -1}\, .
\end{equation}
Using \eqref{eq:dg1}, \eqref{eq:dg2}, and \eqref{eq:dg3}, we see that, for 
$\sigma\in\{+; -\}$,  
\[\int_\C\partial _{\overline{z}}
\varphi^\C (z)(p_z\# p_{\sigma , z}-p_zp_{\sigma , z})\, dz\wedge d\overline{z}\]
converges in $\Sc (h\langle \xi\rangle^{-4}, g)=\Sc (\langle x\rangle^{-1}
\langle \xi\rangle^{-5}, g)$. Thanks to \eqref{eq:caract-pseudo-compact}, 
\begin{equation*}
B_1\ \simeq\ -\sum _{\sigma =\pm}\frac{ikq}{4\pi}\left(\int_\C
\partial _{\overline{z}}\varphi^\C (z)(z-|\xi |^2)^{-1}(z-a_\sigma (x, \xi))^{-1}\, 
dz\wedge d\overline{z}\right)^w e_\sigma\, .
\end{equation*}
We take $b\in\Sc (1, g)$ such that $b\cchi_1=b$. By the previous arguments, 
\begin{align}
b^wB_1&\simeq -\sum _{\sigma =\pm}\left(\int_\C\partial _{\overline{z}}
\varphi^\C (z)b(x, \xi )(z-|\xi |^2)^{-1}(z-a_\sigma (x, \xi))^{-1}\, dz\wedge 
d\overline{z}\right)^w\nonumber\\
&\hspace{5.5cm}\cdot \, ikq(4\pi)^{-1}e^{i\sigma k|x|}\, .\label{eq:int^w}
\end{align}
Now we choose $\varphi $ with a small enough support near some $\lambda\in\Ic$ and 
$b\in\Sc (1, g)$ such that $b(x, \xi)=\cchi _4(x)b_0(\hat{x}, \xi)$, 
$\cchi _4\in\Cc ^\infty(\R^d)$ with $\cchi _4=0$ near $0$ and $\cchi _4=1$ near infinity, 
$\varphi (|\xi |^2)b_0(\hat{x}, \xi)=0=\varphi (|\xi +k\hat{x}|^2)b_0(\hat{x}, \xi)$, 
$b_0=0$ near $\xi \cdot \hat{x}=\pm k/2$, and such that $\varphi (|\xi -k\hat{x}|^2)
b_0(\hat{x}, \xi)$ is nonzero, see Figure \ref{3}.
\begin{figure}
\setlength{\unitlength}{0.6cm}
\begin{center}
\begin{picture}(22,9)(-7,-4.5)
%
\put(-7,0){\vector(1,0){22}} 
\put(0,-0.1){\line(0,1){0.2}} 
\put(3,-0.1){\line(0,1){0.2}}
\put(6,-0.1){\line(0,1){0.2}}
\put(-3,-0.1){\line(0,1){0.2}}
\put(-6,-0.1){\line(0,1){0.2}}
\put(0,-0.6){0}
\put(5.8,-0.6){$k \hat x$}
\put(-6.2,-0.6){$-k \hat x$}
\put(-6,0) {\arc(7,0){-31}}
\put(-6,0) {\arc(7,0){25}}
\put(-6,0) {\arc(7.5,0){-31}}
\put(-6,0) {\arc(7.5,0){25}}
\put(0,0) {\arc(7,0){-31}}
\put(0,0) {\arc(7,0){25}}
\put(0,0) {\arc(7.5,0){-31}}
\put(0,0) {\arc(7.5,0){25}}
\put(6,0) {\arc(7,0){-31}}
\put(6,0) {\arc(7,0){25}}
\put(6,0) {\arc(7.5,0){-31}}
\put(6,0) {\arc(7.5,0){25}}
\put(-3.3,-3.3){\line(0,1){7.2}}
\put(-2.7,-3.3){\line(0,1){7.2}}
\put(3.3,-3.3){\line(0,1){7.2}}
\put(2.7,-3.3){\line(0,1){7.2}}
\put(5,3.5){supp $\varphi(|\cdot|^2)$}
\put(10,3.5){supp $\varphi(|\cdot- k\hat x|^2)$}
\put(-2,3.5){supp $\varphi(|\cdot+ k\hat x|^2)$}
\put(-4.5,-4){$|2 \hat x\cdot \xi +k|< \varepsilon$}
\put(1.5,-4){$|2 \hat x\cdot \xi -k|< \varepsilon$}
\put(13,0){\circle{2}}
\put(10,-1){supp $b_0$}
\end{picture}
\end{center}
\caption{supp $b_0$}\label{3}
\end{figure}
 In the last requirement, we use the 
fact that $\Ic\subset]0;+\infty[$. Note that, on the support of $b_0(\hat{x}, \xi)$ and for 
$|x|$ large enough, $b_1(x, \xi ):=|\xi |^2-|\xi +\sigma k\hat{x}|^2$ does not vanish, for 
$\sigma\in\{+;-\}$. Thus, 
\begin{align*}
(z-|\xi |^2)^{-1}(z-a_\sigma (x, \xi))^{-1}=(b_1(x, \xi ))^{-1}\bigl((z-|\xi |^2)^{-1}-
(z-a_\sigma (x, \xi))^{-1}\bigr)\, ,
\end{align*}
in this region. Inserting this in \eqref{eq:int^w} and using the support properties of $b$ and $\varphi$, 
\begin{align*}
b^wB_1&\simeq -kq2^{-1}\sum _{\sigma =\pm}\Bigl(b(x, \xi )(b_1(x, \xi ))^{-1}
\bigl(\varphi (|\xi |^2)-\varphi (|\xi +\sigma
k\hat{x}|^2)\bigr) \Bigr)^we^{i\sigma k|x|}\\ 
&\simeq kq2^{-1}\Bigl(b(x, \xi )(b_1(x, \xi ))^{-1}
\varphi (|\xi - k\hat{x}|^2)\Bigr)^we^{-ik|x|}\, .
\end{align*}
Setting $B=b^w$ and $B'=e^{ik|x|}$, $BB_1B'\simeq c^w$
with an explicit $c\in\Sc (1, g)$ that does not tend to $0$ at
infinity. By \eqref{eq:caract-pseudo-compact}, neither 
$c^w$ nor $BB_1B'$ is compact. \qed
\begin{rem}\label{rq:autre-op-conj}
As alternative to $A_1$, it is natural to try $\hat{A}_1=(\tau
(\xi)\rho (x)x\cdot\xi )^w$,  where $1-\rho \in \Cc ^\infty_c(\R^d)$
and $\tau \in S(1, g)$ satisfies $\tau (\xi )=1$ if  $|\xi
|^2\in\Ic$. But Proposition 5.3 holds true with $A_1$ replaced by
$\hat{A}_1$.  
\end{rem}
Let us sketch a justification of Remark~\ref{rq:autre-op-conj}. One can
verify that  $\varphi (H_0)\in\Cc^2(\hat{A}_1)$. We follow the
proof of Lemma~\ref{l:non-compact} and arrive at
\eqref{eq:B1-equiv-loin-0} where $\chi _1W_1$ is replaced by:  
$$-\sum _{\sigma =\pm}\frac{ikq}{4\pi}\int_\C \partial _{\bar{z}}
\varphi ^\C(z)(z-H_0)^{-1}b_{\tau , \sigma}^w\, dz\wedge d\bar{z}\, e_\sigma$$
with, for $\chi _2=0$ near $0$ and $\chi _2\chi _1=\chi _1$, 
$$b_{\tau , \sigma}=\chi _2(x)\bigl(\tau (\xi)\hat{x}\cdot\xi -\tau (\xi -\sigma k\hat{x})\hat{x}\cdot
(\xi-\sigma k\hat{x})\bigr)\, .$$
Since the $b_{\tau , \sigma}$ do not depend on $z$, we can estimate
$p_z\# b_{\tau , \sigma}\#  p_{\sigma , z}-b_{\tau , \sigma}p_zp_{\sigma , z}$ in a similar 
way as in \eqref{eq:esti-z} and get \eqref{eq:int^w} with  $b$
replaced by $bb_{\tau , \sigma}$. Following the last lines, we find
that $BB_1B'\simeq  (b_{\tau , -}c)^w$, $b_{\tau , -}c\in\Sc (1, g)$, and $b_{\tau , -}c$ does 
not tend to zero at infinity. We arrive at the same conclusion as in 
Lemma~\ref{l:non-compact}.

\appendix 
\renewcommand{\theequation}{\thesection .\arabic{equation}}

\section{Oscillating terms.} 
\label{app:s:oscillation}
\setcounter{equation}{0}

In our study of Schr\"odinger operator with a perturbed Wigner-Von Neumann potential 
(see Section~\ref{s:wigner}), we need a good understanding of operator compositions like 
$a^w\cchi _1W_1$, where $a\in\Sc (m,g)$, $g$ and $m$ given by \eqref{eq:metric} and  
\eqref{eq:weight}, $W_1$ given by \eqref{eq:def-W_1}, and where $\cchi _1\in\Cc ^\infty (\R^d)$ 
such that $\cchi _1=0$ near $0$ and $\cchi _1=1$ near infinity. More precisely, we are looking 
for an explicit pseudodifferential operator $A$ such that $a^w\cchi _1W_1=A+b_1^wB_1+B_2b_2^w$, 
with bounded operators $B_1, B_2$ and symbols $b_1, b_2\in\Sc (m\langle x\rangle ^{-1}\langle \xi\rangle ^{-1},g_0)$ 
($g_0$ given in \eqref{eq:metric}). Although $a\in\Sc (m,g_0)$ and $\cchi _1W_1\in\Sc 
(\langle x\rangle ^{-1},g_0)$, the symbolic calculus associated to $g_0$ is not well 
suited for our analysis, in particular to guarantee $b_1, b_2\in\Sc (m\langle x\rangle ^{-1}\langle 
\xi\rangle ^{-1},g_0)$. It is better to work with $g$ with the drawback that $W_1$ does not 
belong the corresponding calculus. Taking into account the special form of $W_1$, we provide the previous 
decomposition with $b_1, b_2\in\Sc (m\langle x\rangle ^{-1}\langle \xi\rangle ^{-1},g)$, using 
standard arguments of pseudodifferential calculus. In Appendix~\ref{app:s:dim1}, we give a 
simpler result in dimension $d=1$ that essentially follows from facts used in \cite{fh}.

For $m$ of the form \eqref{eq:weight}, we denote by $\Sc (m\langle x\rangle ^{-\infty}, g)$ 
the intersection of all classes $\Sc (m\langle x\rangle ^{k}, g)$ for $k\in\Z$. We denote by 
$\Sc (-\infty ,g)$ the intersection of all classes $\Sc (m, g)$ with $m$ satisfying 
\eqref{eq:weight}. It suffices to study $a^we_{\pm}$ where 
$e_{\pm}(x)=\cchi _1(x)e^{\pm ik|x|}$. To this end, we shall use the oscillatory 
integrals defined in Theorem 7.8.2, p.\ 237, in \cite{h1}, which actually works 
for symbols in the classes $\Sc (m,g)$ we consider here. These oscillatory 
integrals can also be viewed as tempered distributions. Note that usual operations 
on integrals (like integration by parts or change of variable) are valid for 
oscillatory integrals. 

\begin{proposition}\label{p:a^w-oscillant}
Let $a\in\Sc (m, g)$ with $m$ and $g$ given by \eqref{eq:weight} and \eqref{eq:metric}. 
Let $e_{\pm}$ be the functions defined just above. Then there exist symbols $a_\pm\in
\Sc (m, g)$, $b_\pm\in \Sc (m\langle x\rangle ^{-\infty}, g)$, and $c_\pm
\in\Sc (mh, g)$ (with $h$ defined in \eqref{eq:gain}), such that $e_{\pm}a^w=a_{\pm}^we_{\pm}
+e_{\pm}b_\pm^w+c_\pm^we_{\pm}$ and such that $a_{\pm}(x, \xi )=a(x, \xi \mp k|x|^{-1}x)$, 
if $\cchi _1(x)\neq 0$. 
\end{proposition} 
\proof  Let $\cchi _2, \check{\cchi }_2\in\Cc ^\infty (\R^d)$ such that $\cchi _2=0$ and 
$\check{\cchi }_2=0$ near $0$, $\check{\cchi }_2\cchi _1=\cchi _1$, and $\check{\cchi }_2
(1-\cchi _2)=0$. 
Notice that $\cchi _2\cchi _1=\cchi _1$. We write $e_{\pm}a^w=e_{\pm}a^w(\cchi _2+1-\cchi_2)=
e_{\pm}a^w\cchi _2e^{\mp ik|x|}e_{\pm}+e_{\pm}\check{\cchi }_2a^w(1-\cchi_2)$ and arrive at 
\begin{equation}\label{eq:epma^w}
e_{\pm}a^w\ =\  e_{\pm}a^w\cchi _2e^{\mp ik|x|}
e_{\pm}\, +\, e_{\pm}b^w
\end{equation}
where $b:=\check{\cchi }_2\# a\# (1-\cchi _2)\in\Sc (-\infty  ,g)$, since 
$\check{\cchi }_2(1-\cchi _2)=0$. For $f\in\Sr (\R^d)$, the Schwartz
space on $\R^d$, using an  oscillatory integral in the $\xi$ variable, 
\begin{align*}
f_1(x)&:=(e_{\pm}a^w\cchi _2e^{\mp ik|x|}f) (x)\\
&=(2\pi )^{-d}\int e^{i\langle x-y, \xi\rangle}a((x+y)/2;\xi )\cchi
_1(x)e^{\pm ik|x|}\cdot \ \cchi _2(y)
e^{\mp ik|y|}f(y)\, dyd\xi\\ 
&=(2\pi )^{-d}\int e^{i\langle x-y, \xi\rangle\pm
  ik(|x|-|y|)}a((x+y)/2;\xi )\cdot \  \cchi _1(x)\cchi _2(y)f(y)\,
dyd\xi \, .\\ 
\end{align*}
We take $\varepsilon\in]0;1/4[$ and $\tau\in\Cc^\infty_c (\R)$ such
that $\tau (t)=1$ if $|t|\leq  1-4\varepsilon$ and $\tau (t)=0$ if
$|t|\geq 1-2\varepsilon$. We insert $\tau (|x-y|\langle
x\rangle^{-1})+ 1-\tau (|x-y|\langle x\rangle^{-1})$ into the previous
expression of $f_1$ and call $f_2$ (resp.\   $f_3$) the integral
containing $\tau (|x-y|\langle x\rangle^{-1})$ (resp.\  $1-\tau (|x-y| 
\langle x\rangle^{-1})$). On the support of $\cchi _1(x)\cchi _2(y)\tau (|x-y|\langle x\rangle^{-1})$, 
$|x-y|\leq (1-2\varepsilon )\langle x\rangle$. We can choose the support of $\cchi _1$ such that, 
on the support of $\cchi _1(x)\cchi _2(y)\tau (|x-y|\langle x\rangle^{-1})$, $|x-y|\leq (1-\varepsilon )
|x|$. In particular, on this support, $0$ does not belong the segment $[x;y]$ and, for all 
$t\in [0;1]$, 
\begin{equation}\label{eq:def-u}
 u(t; x, y):= |tx+(1-t)y|\geq |x|-(1-t)|y-x|\geq \varepsilon |x|\, .
\end{equation}
For $x\neq y$, $(L_{x, y, D_\xi}-1)e^{i\langle x-y, \xi\rangle\pm ik(|x|-|y|)}=0$ for 
$L_{x, y, D_\xi}=|x-y|^{-2}(x-y)\cdot D_\xi$. Thus, by integration by parts, for all 
$p\in\N$,
\begin{align}
f_3(x)&=(2\pi )^{-d}\int e^{i\langle x-y, \xi\rangle\pm ik(|x|-|y|)}
\cchi _1(x)\cchi _2(y)(1-\tau (|x-y|
\langle x\rangle^{-1}))\nonumber\\
&\hspace{4.3cm}\cdot \ \bigl(L_{x, y, D_\xi}^\ast\bigr)^p\bigl(a((x+y)/2;\xi )
\bigr)f(y)\, dyd\xi \nonumber\\
\label{eq:f_3}&=(b_3^wf)(x)\, ,
\end{align}
with $b_3\in\Sc (-\infty  ,g)$ (cf., (8.1.8) in \cite{h3}). 
\begin{lemma}\label{l:prop-u}
Take $x, y\in\R^d$ such that $0$ does not belong the segment $[x;y]$. Then,
\begin{equation}\label{eq:|x|-|y|}
|x|-|y|\ =\ \langle v(1/2; x, y)+r(x, y), x-y\rangle
\end{equation}
where $v(t; x, y)=(tx+(1-t)y)/|tx+(1-t)y|$ for $t\in [0;1]$, and where 
\begin{equation}\label{eq:r}
r(x, y):=\int \bigl((1-t)\un _{[1/2;1]}(t)-t\un
_{[0;1/2]}(t)\bigr)\partial _tv (t; x,y )\, dt 
\end{equation}
satisfies $|r(x, y)|\leq 2$. 
\end{lemma}
\proof It suffices to use the Taylor expansion with integral rest for the function $u(\cdot ; x, y)$ 
defined in \eqref{eq:def-u} between $0$ and $1/2$ and between $1/2$ and $1$. \qed

By Lemma~\ref{l:prop-u}, we can rewrite $f_2(x)$ as 
\begin{align*}
f_2(x)&=(2\pi )^{-d}\int e^{i\langle x-y, \xi\pm k(v(1/2; x, y)+r(x, y))\rangle}
\cchi _1(x)\cchi _2(y)\tau (|x-y|
\langle x\rangle^{-1})\\
&\hspace*{5cm}\cdot \ a((x+y)/2;\xi )f(y)\, dyd\xi \\
&=(2\pi )^{-d}\int e^{i\langle x-y, \eta\rangle}\cchi _1(x)\cchi _2(y)\tau (|x-y|
\langle x\rangle^{-1})\\
&\hspace*{3cm}\cdot \ a((x+y)/2;\eta \mp  k(v(1/2; x, y)+r(x, y)))f(y)\, dyd\eta \, ,
\end{align*}
after the change of variable $\eta =\xi\pm k(v(1/2; x, y)+r(x, y))$. 
Now we use a Taylor expansion of $a$ with integral rest in the $\xi$ variable: 
\begin{align*}
&\hspace*{-0.5cm}a\big((x+y)/2;\eta \mp  k(v(1/2; x, y)+r(x, y))\big)
= a((x+y)/2;\eta \mp  kv(1/2; x, y))\\ 
&\hspace{1cm}+\int _0^1dt\, \langle\nabla _\xi a ((x+y)/2;\eta \mp
k(v(1/2; x, y)+tr(x, y))),  
kr(x, y)\rangle \, .
\end{align*}
According to this decomposition, we split $f_2(x)$ in $f_4(x)+f_5(x)$. Thanks to the bound 
\eqref{eq:def-u} for $t=1/2$, we can find $\cchi _3\in\Cc ^\infty (\R^d)$ such that 
$\cchi _3=0$ near $0$ and 
\[\cchi _1(x)\cchi _2(y)\tau (|x-y|\langle x\rangle^{-1})(1-\cchi _3((x+y)/2))=0\, .\]
Setting $a_\pm (x, \eta )=\cchi _3(x)a(x, \eta \mp k\hat{x})$, we obtain that 
\begin{align}
f_2(x)&=(2\pi )^{-d}\int e^{i\langle x-y, \eta\rangle}\cchi _1(x)\cchi
_2(y)\tau (|x-y| 
\langle x\rangle^{-1})\nonumber\\
&\hspace{5cm}\cdot \ a_\pm((x+y)/2;\eta )f(y)\, dyd\xi +f_5(x)\nonumber\\
&=\cchi _1(x)(a_\pm^w\cchi _2f)(x)+f_5(x)
\, =\, (a_\pm^wf)(x)+(b_2^wf)(x)+f_5(x)\, , \label{eq:f_2}
\end{align}
with $b_2\in\Sc (m\langle x\rangle ^{-\infty}, g)$. Since $||\eta +k\hat{x}|-|\eta||\leq k$, 
for all $x\in\R^d\setminus\{0\}$ and all $\eta\in\R^d$, 
a direct computation shows that $a_\pm\in\Sc (m,g)$. \\
Now we study $f_5$. Given a vector $v\in\R^d$, let $A(v)=I-\langle v, \cdot\rangle v$
(where $I$ denotes the identity on $\R^d$). If $0$ does not belong to the segment 
$[x; y]$ in $\R^d$, $\partial _tv(t; x, y)=(u(t; x, y))^{-1}A(v(t; x, y))\cdot (x-y)$, where 
$v(t; x, y)$ (resp.\  $u(t; x, y)$) is defined in Lemma~\ref{l:prop-u}
(resp.\  \eqref{eq:def-u}). Defining $\kappa (s):=(1-s)\un _{[1/2;1]}(s)-s\un _{[0;1/2](s)}$, 
\begin{align*}
f_5(x)&=(2\pi )^{-d}\int e^{i\langle x-y, \eta\rangle}\cchi _1(x)\cchi _2(y)\tau (|x-y|
\langle x\rangle^{-1})\int_0^1dt\\
&\hspace{4cm}\cdot \, \bigl\langle \nabla _\xi a((x+y)/2;\eta \mp k(v(1/2; x, y)+t
r(x, y)))\, ,\\
&\hspace{2cm}\cdot \, k\int_0^1ds\,  \kappa (s)(u(s; x, y))^{-1}A(v(s;
x, y))\cdot (x-y) 
\bigr\rangle\cdot \,  f(y)\, dyd\xi \, , 
\end{align*}
by \eqref{eq:r}. 
Denoting by $A(v)^T$ the transposed of the linear map $A(v)$ and setting $\eta _t=\eta \mp 
k(v(1/2; x, y)+tr(x, y))$, 
\begin{align*}
&\bigl\langle \nabla _\xi a((x+y)/2;\eta _t )\, ,\,  A(v(s; x, y))\cdot (x-y)
\bigr\rangle \\
&\hspace*{5cm}=\bigl\langle A(v(s; x, y))^T\nabla _\xi a((x+y)/2;\eta _t)\, ,\,  (x-y)\bigr\rangle\, . 
\end{align*}
Integrating by parts in the $\eta$ variable, 
\begin{align*}
f_5(x)&=(2\pi )^{-d}\int e^{i\langle x-y, \eta\rangle}\cchi _1(x)\cchi _2(y)\tau (|x-y|
\langle x\rangle^{-1})\int_0^1dt\int_0^1ds\\
&\hspace{3.5cm}\cdot \, \bigl(i\langle A(v(s; x, y))\nabla _\xi , \nabla _\xi\rangle a\bigr)
((x+y)/2;\eta _t )\\
&\hspace{3.5cm}\cdot \, k\kappa (s)(u(s; x, y))^{-1}\, f(y)\, dyd\xi \, . 
\end{align*}
Writing $f(y)=(2\pi )^{-d}\int e^{i\langle y, \xi\rangle}(\Fc f)(\xi )d\xi$,  
$\Fc f$ being the Fourier transform of $f$, 
\begin{equation}\label{eq:f_5}
f_5(x)\ =\ (2\pi )^{-d}\int e^{i\langle x, \xi\rangle}c_0(x, \xi)(\Fc f)(\xi )d\xi \, =\, 
(\op c_0f)(x)
\end{equation}
where $c_0$ is defined by the oscillatory integral (in the $\eta$ variable) 
\begin{align}
\label{eq:c_0}
c_0(x, \xi)\ &=\ \int e^{i\langle x-y, \eta -\xi\rangle}\rho (x, y; \eta )\, dyd\eta
\hspace{0.5cm}\mbox{with}\\
\rho (x, y; \eta )&=\cchi _1(x)\cchi _2(y)\tau (|x-y|
\langle x\rangle^{-1})\int_0^1dt\int_0^1ds\, (u(s; x, y))^{-1}\nonumber\\
&\hspace*{1cm}\label{eq:def-rho}\cdot \, ik\kappa (s)\bigl(\langle A(v(s; x, y))
\nabla _\xi , \nabla _\xi\rangle a\bigr)((x+y)/2;\eta _t )\, . 
\end{align}
Now we inset in \eqref{eq:c_0} $\tau (|\eta -\xi |\langle \xi\rangle^{-1})+1-\tau (|\eta -\xi |
\langle \xi\rangle^{-1})$ and split $c_0$ into $c_1+c_2$. In particular, 
\begin{align*}
 c_2(x, \xi)&= \int e^{i\langle x-y, \eta -\xi\rangle}\bigl(1-\tau (|\eta -\xi |
\langle \xi\rangle^{-1}\bigr)\rho (x, y; \eta )\, dyd\eta\\
&=\int e^{i\langle x-y, \eta -\xi\rangle}\bigl(1-\tau (|\eta -\xi | 
\langle \xi\rangle^{-1}\bigr)\bigl(L_{\xi , \eta ,
  D_y}^\ast\bigr)^p(\rho (x, y; \eta ))\, dyd\eta 
\end{align*}
for all $p\in\N$. By direct computations, we see that $c_2\in\Sc
(-\infty , g)$ and  $c_1\in\Sc (mh,g)$. Since for any symbol $r$,
there exists a symbol $s$ in the same class  such that $\op r=s^w$,
the equations \eqref{eq:epma^w}, \eqref{eq:f_3}, \eqref{eq:f_2},  
and \eqref{eq:f_5}, yield the desired result. \qed

\section{Functional calculus for pseudodifferential operators. } 
\label{app:s:calcul-fonct-pseudo}
\setcounter{equation}{0}

Here we present a result on the functional calculus for pseudodifferential operators associated 
to the metric $g$ in \eqref{eq:metric}. This result is probably not new but we did not find a proof 
in the literature. It follows quite directly from arguments in \cite{bo} (see also \cite{l}). 
However we sketch the proof for completeness. We use notions and results from 
Subsection~\ref{ss:psi-do}. 

Recall that, for $\rho\in\R$, we denote by $\Sc^\rho$ the set of smooth functions $\varphi$ on 
$\R$ such that $\sup _{t\in\R}\langle t\rangle ^{k-\rho}|\partial _t^k\varphi (t)|<\infty$. 
If we take a real symbol $a\in \Sc (m,g)$, then the operator $a^w$ is
self-adjoint on the  
domain $\Dc (a^w)=\{u\in\rL^2(\R^d_x);a^wu\in\rL^2(\R^d_x)\}$. In particular, the operator 
$\varphi (a^w)$ is well defined by the functional calculus if $\varphi$ is a borelean 
function on $\R$. We assume that $m\geq 1$. A real symbol $a\in \Sc (m,g)$ is said elliptic 
if $(i-a)^{-1}$ belongs to $\Sc (m^{-1},g)$. 
\begin{theorem}\label{t:func}
Let $m\geq 1$ and $a\in \Sc(m, g)$ be real and elliptic. Take $\varphi\in \Sc ^\rho$. Then 
$\varphi (a)\in \Sc(m^\rho , g)$ and there is $b\in \Sc(hm^\rho , g)$ such that 
\begin{equation}\label{eq:fonct-pseudo}
\varphi\big(a^w(x,D)\big)= \big(\varphi(a)\big)^w(x,D) + b^w(x,D).
\end{equation}
\end{theorem}
\proof Let $\rho '\in\R$, $\varphi\in \Sc^{\rho '}$, and $k\in\N$ large enough such that $2k>\rho '$.  
Then $\psi (t):=\varphi (t)(1+t^2)^{-k}$ belongs to $\Sc^{\rho '-2k}$ with $\rho '-2k<0$. 
If the result is valid for $\rho <0$, then there exists $b\in\Sc (hm^{\rho '-2k} , g)$ such that 
\begin{align*}
 \varphi (a^w)&=\psi (a^w)(1+(a^w)^2)^k=((\psi (a))^w+b^w)(1+(a^w)^2)^k\\
&=\bigl(\psi (a)(1+a^2)^k\bigr)^w+c^w+(b\# (1+a^2)^k)^w=(\varphi (a))^w+d^w
\end{align*}
with $c, d\in \Sc (hm^{\rho '} , g)$, by the composition properties. So it suffices to 
prove the result for $\rho<0$. Since we can write any function $\varphi\in \Sc^\rho$, 
with $\rho<0$, as $\varphi _1\varphi _2$ with $\varphi _1\in \Sc^\delta$ ($-1\leq\delta <0$) and 
$\varphi _1\in \Sc^{[\rho]+1}$ (where $[\rho]$ denotes the integer part of $\rho$) and 
use the previous composition properties, we see by induction that it suffices to 
establish the result for $-1\leq\rho<0$.\\
Let $z\in\C\setminus \R$. Using that $(z-a)^{-1}=(i-a)^{-1}(1+(z-i)(z-a)^{-1})$, 
we observe that $|(z-a)^{-1}|\leq m^{-1}\langle z\rangle |\im
(z)|^{-1}$. Thus, for all  
$\ell\in\N$, there exists $C_\ell >0$ and $N_\ell \in\N$ such that, 
for all $z\in \C\setminus \R$, $\|(z-a)^{-1}\|_{\ell ,\Sc (m^{-1},g)}\leq C_\ell 
\langle z\rangle^{N_\ell +1}|\im (z)|^{-N_\ell -1}$. Define 
$q_z:=(z-a)^{-1}\# (z-a)-1\in \Sc (h,g)$.  By an explicit formula given
in \cite{bo2} (first formula on page III-4), $q_z$ only depends on the
derivatives of $(z-a)$, which are independent of $z$. Thus one can
find, for all $\ell$, $C_\ell '>0$ and  $N_\ell '\in\N$ such that 
\begin{equation}\label{eq:bound-q_z}
\forall z\in \C\setminus \R\, ,\hspace{.3cm}\|q_z\|_{\ell , \Sc (h,g)}\leq C_\ell '\langle z\rangle^{N_\ell '+1}|\im (z)|^{-N_\ell '-1}\, .
\end{equation}
According to \cite{bo}, one can prove from the boundedness of commutators of 
$(z-a^w)^{-1}$ with appropriate pseudodifferential operators that this resolvent $(z-a^w)^{-1}$
is equal to $r_z^w$, where the symbol $r_z$ belongs to $\Sc (m^{-1},g)$. Furthermore, 
the system $\|\cdot \|_{\ell , \Sc (m^,g)}$, $\ell\in\N$, of semi-norms 
is equivalent to another one based on the previous commutators. 
Using this, there exist, for all $\ell$, $C_\ell ''>0$ and 
$N_\ell ''\in\N$ such that
\begin{equation}\label{eq:bound-r_z}
\forall z\in \C\setminus \R\, ,\hspace{.3cm}\|r_z\|_{\ell , \Sc (m^{-1},g)}\leq C_\ell ''\langle z\rangle^{N_\ell ''+1}|\im (z)|^{-N_\ell ''-1}\, .
\end{equation}
Using \eqref{eq:bound-q_z} and \eqref{eq:bound-r_z}, we can find, for all $\ell$, $C_\ell '''>0$ and 
$N_\ell '''\in\N$ such that
\begin{equation}\label{eq:bound-q_zdièser_z}
\forall z\in \C\setminus \R\, ,\hspace{.3cm}\|q_z\# r_z\|_{\ell , \Sc (hm^{-1},g)}\leq C_\ell '''\langle z\rangle^{N_\ell '''+1}
|\im (z)|^{-N_\ell ''-1}\, .
\end{equation}
Now we take $\varphi\in\Sc^\rho$ with $-1\leq\rho<0$ and consider some almost analytic 
extension $\varphi ^\C$ (like in Proposition~\ref{p:dg}). Thanks to 
\eqref{eq:bound-q_zdièser_z}, \eqref{eq:dg1}, \eqref{eq:dg2}, and $\rho <0$, 
\[b:=\frac{i}{2\pi}\int _\C\partial_{\bar{z}}\varphi ^\C(z)q_z\# r_z\, dz\wedge 
d\overline{z}\]
converges in $\Sc (hm^{-1}g)$. According to the definition of $q_z$, 
$((z-a)^{-1})^w(z-a^w) = {\rm Id} 
+ q_z^w$, thus $((z-a)^{-1})^w = (z-a^w)^{-1} + (q_z\# r_z)^w$. Using  
Helffer-Sj\"ostrand formula \eqref{eq:int}, $(\varphi (a))^w=\varphi (a^w)+b^w$ with 
$b\in\Sc (hm^{-1}g)\subset \Sc (hm^{\rho}g)$, since $-1\leq\rho$. \qed

\section{An interpolation's argument. } 
\label{app:s:interpolation}
\setcounter{equation}{0}  

By pseudodifferential calculus, $A_1^2\langle D_x\rangle^{-2}\langle
x\rangle^{-2}$ extends  
to a bounded operator on $\rL^2(\R _x^d)$. What about $\langle A_1\rangle^r
\langle D_x\rangle^{-r}\langle x\rangle^{-r}$ with $r>0$? The same
argument is not clear since $A_1$ is not elliptic. Indeed its symbols
$(x,\xi )\donne x\cdot\xi$ can vanish when $\xi\neq 0$. Using
interpolation, we show  
\begin{lemma}\label{l:interpolation}
For real $r\geq 0$, 
$\langle A_1\rangle^r\langle D_x\rangle^{-r}
\langle x\rangle^{-r}$ extends to a bounded operator on $\rL^2(\R _x^d)$. 
\end{lemma}
We refer to \cite{ms}[Lemma 7.1] for an alternative proof and
historical remarks. 

\proof[Proof of Lemma \ref{l:interpolation}] We prove that, for $r\geq 0$, 
\begin{equation}\label{eq:norm_r}
\exists C_r>0;\, \forall f\in\rL^2(\R^d_x),\, \|\langle A_1\rangle^r f\|\leq 
C_r \|\langle x\rangle^{r} \langle D_x\rangle^{r} f\|\, .
\end{equation}
For $r\in\N$, $(A_1+i)^r\langle D_x\rangle^{-r}
\langle x\rangle^{-r}$ extends to a bounded operator by the
pseudodifferential calculus  
with the metric $g$ in \eqref{eq:metric}. Since $\langle A_1\rangle^r
(A_1+i)^{-r}$ is  
bounded, \eqref{eq:norm_r} is satisfied when $r\in\N$. For $t,t'\geq 0$, let 
$\rH^{t'}_t:= \{f\in\rL^2(\R^d_x);\|\langle x\rangle^{t}\langle
D_x\rangle^{t'}f\|<\infty\}$. 
Now, using \cite{be}, we infer that the space $\rH^r_r$ is also the complex 
interpolated space $[\rH^0_0, \rH^m_m]_{r/m}$, where $m\geq r$. To be precise, 
use \cite{be}[(1.7)] and notice that 
$H(m,g)=\rH^r_r$, where $g$ is given as in \eqref{eq:metric} and 
$m(x,\xi)=\langle x\rangle^{r} \langle \xi\rangle^{r}$, by \cite{be}[Theorem 3.7]. We deduce that 
\eqref{eq:norm_r} is true for all $r\geq 0$ by the Riesz-Thorin Theorem.\qed

\section{A simpler argument in dimension $d=1$. } 
\label{app:s:dim1}
\setcounter{equation}{0}  

Here we present a more elementary proof of Lemma~\ref{l:non-compact} in dimension $d=1$. 
It relies on the following 
\begin{lemma}\label{l:f(p)-oscillant-dim-1}
For $z\not\in\R$, as bounded operators on $\rL ^2(\R_x)$, 
\begin{equation}\label{eq:interversion-dim-1}
(z-D_x^2)^{-1}e^{\pm ikx}\ =\ e^{\pm ikx}\bigl(z-(D_x\pm k)^2\bigr)^{-1}\, .
\end{equation}
\end{lemma}
\proof As differential operators, $D_xe^{\pm ikx}=e^{\pm ikx}(D_x\pm k)$. Thus, 
on $\rH^2 (\R_x)$, $(D_x^2-z)e^{\pm ikx}=e^{\pm ikx}((D_x\pm k)^2-z)$. Multiplying 
on the left and on the right by the convenient resolvent, we get the result. \qed

We first follow the general proof until formula~\eqref{eq:B1-equiv}. By 
\eqref{eq:interversion-dim-1}, 
\[B_1\ \simeq\ -\sum_{\sigma =\pm}\frac{iqk}{4\pi}\int_\C\partial _{\overline{z}}\varphi^\C (z)
(z-D_x^2)^{-1}(z-(D_x+\sigma k)^2)^{-1}\, dz\wedge d\overline{z}\, e^{i\sigma kx}\, .\]
Choosing the support of $\varphi$ small enough, we can find $\theta\in\Cc^\infty _c(\R)$ such 
that $\xi\donne\theta (\xi ^2)$ vanishes near $-k/2$ and $k/2$, 
$\varphi (\xi^2)\theta (\xi^2)=0=\varphi ((\xi +k)^2)\theta (\xi^2)$, for all $\xi\in\R$, 
and such that the function $\xi\donne \varphi ((\xi -k)^2)\theta (\xi^2)$ is nonzero 
(using that $\Ic\subset ]0;+\infty [$). Set $B=\theta (D_x^2)$. 
Since $\xi^2-(\xi +k)^2$ and $\xi^2-(\xi -k)^2$ do not vanish on the support of 
$\theta (\xi^2)$, 
\begin{align*}
BB_1&\simeq -\sum_{\sigma =\pm}\frac{iqk}{4\pi}\theta
(D_x^2)(D_x^2-(D_x +k\sigma )^2)^{-1} \int_\C\partial
_{\overline{z}}\varphi^\C (z)(z-D_x^2)^{-1}\\ 
 &\hspace{3cm}(D_x^2-(D_x +k\sigma )^2)(z-(D_x+\sigma k)^2)^{-1}\, dz\wedge 
d\overline{z}\, e^{i\sigma kx}\, .
\end{align*}
By the resolvent formula and \eqref{eq:int}, 
\[BB_1\ \simeq\ -\sum_{\sigma =\pm}\frac{qk}{2}\theta (D_x^2)(D_x^2-(D_x +k\sigma )^2)^{-1}
\bigl(\varphi (D_x^2)-\varphi ((D_x+k\sigma )^2)\bigr)\, e^{i\sigma kx}\, .\]
Using the support properties of $\theta$, we obtain 
\[BB_1\ \simeq\ 2^{-1}qk\theta (D_x^2)(D_x^2-(D_x -k)^2)^{-1}\varphi ((D_x-k)^2)
\, e^{-ikx}\, .\]
Denoting by $B'$ the multiplication operator by $e^{ikx}$, $BB_1B'$
is, modulo some compact  operator, a self-adjoint Fourier
multiplier. The spectrum of the latter is given by the essential range
of the function $\xi\donne 2^{-1}qk\theta (\xi^2)(\xi^2-(\xi
-k)^2)^{-1}  \varphi ((\xi-k)^2)$. Since this
function is non constant and continuous,  the spectrum contains an
interval and the corresponding operator cannot be compact. Thus
$BB_1B'$ is not compact. This finishes the proof of
Lemma~\ref{l:non-compact} in dimension $d=1$.  \qed

  
%

\begin{thebibliography}{xxxxxx}  


%
\bibitem[Ag]{a}S.\ Agmon: {\em Lower bounds for solutions of
  Schr\"odinger equations.} J.\ Analyse Math.\  23, 1970, 1-25.
%
\bibitem[ABG]{abg}W.O.\ Amrein, A.\ Boutet de Monvel, V.\ Georgescu: {\em
 $C_0$-groups, commutator methods and spectral theory of $N$-body   
hamiltonians.}, Birkh\"auser 1996.   
%
\bibitem[BFS]{bfs}V.\ Bach, J.\ Fr\"ohlich, I.M.\ Sigal: 
{\em Quantum electrodynamics of confined non-relativistic particles}.
Adv. in Math. {\bf 137}, 299--395, 1998.
%
%
%
\bibitem[Be]{be} R.\ Beals: {\em Weighted distribution spaces and pseudodifferential 
operators}, J.\ Analyse Math.\ 39 (1981), 131--187. 
%
\bibitem[BD]{bd}M.\ Ben-Artzi, A.\ Devinatz: {\em Spectral and
  scattering theory for the adiabatic oscillator and related
  potentials.} J. Maths. Phys. 111 (1979), p.\ 594-607. 
%
\bibitem[BCHM]{bchm}J-F.\ Bony, R. Carles, D. H\"afner, L. Michel: {\em Scattering theory 
for the Schr\"odinger equation with repulsive potential.} J. Math. Pures Appl. 84, no. 5, 509-579, 
2005.
%
\bibitem[Bo1]{bo}J-M.\ Bony: {\em Caract\'erisation des op\'erateurs pseudo-diff\'erentiels.} 
S\'eminaire \'EDP (1996-1997), expos\'e num\'ero XXIII.
%
\bibitem[Bo2]{bo2}J-M.\ Bony: {\em Sur l'inégalit\'e de Fefferman-Phong.} 
S\'eminaire \'EDP (1998-1999), expos\'e num\'ero III.
%
\bibitem[BC]{bc}J-M.\ Bony, J-Y. Chemin: {\em Espaces fonctionnels associ\'es au calcul de 
Weyl-H\"ormander.} Bulletin de la S.M.F., tome 122, num\'ero 1 (1994), p.\ 77-118.
%
\bibitem[BG]{bg} N. Boussaid, S. Gol\'enia: {\em Limiting absorption principle for some 
long range pertubations of Dirac systems at threshold energies.}   
Comm.\ Math.\ Phys.\ 299, 677-708 (2010). 
%
\bibitem[Bu]{b}N.\ Burq: {\em Semiclassical estimates for the resolvent   
in non trapping geometries.} Int.\ Math.\ Res.\ Notices 2002, no 5, 221--241.  
%
\bibitem[CJ]{cj}F.\ Castella, Th.\ Jecko: {\em Besov estimates in the high-frequency 
Helmholtz equation, for a non-trapping and $C^2$ potential.} J. Diff. Eq. Vol. 228, N.2, p.\ 440-485 (2006).  
%
\bibitem[Ca]{c} L.\ Cattaneo: {\em Mourre's inequality and Embedded bounded
  states}, Bull.\ Sci Math.\ {\bf 129}, Issue 7, 591--614, 2005. 
%
\bibitem[CGH]{cgh} L.\ Cattaneo, G.\ M.\ Graf, W.\ Hunziker:
{\em A general resonance theory based on Mourre's inequality}, Ann. H. Poincar\'e 7 (2006), 
3, p.\ 583-601.
%
\bibitem[CHM]{chm}J.\ Cruz-Sampedro, I.\ Herbst, R.\ Martinez-Avenda\~{n}o: {\em   
Pertubations of the Wigner-Von Neumann potential leaving the embedded eigenvalue fixed.} 
Ann. H. Poincar\'e 3 (2002) 331-345.   

%
\bibitem[CFKS]{cfks}H.L.\ Cycon, R.G.\ Froese, W.\ Kirsch, B.\ Simon: {\em Schr\"odinger 
operators with applications to quantum mechanics and   global geometry.} Springer Verlag 1987.   
%


%
\bibitem[DG]{dg}J.\ Derezi\'nski, C.\ G\'erard: {\em Scattering theory of
  classical and quantum N-particle systems.} Springer-Verlag 1997.   
%
\bibitem[DJ]{dj}J.\ Derezi\'nski, V.\ Jak\v{s}i\'c: {\em Spectral theory of Pauli-Fierz 
operators.} J.\ Funct.\ Anal.\ {\bf 180}, no 2, p.\ 243-327, 2001.  
%
\bibitem[DMR]{dmr}A.\ Devinatz, R.\ Moeckel, P.\ Rejto: {\em A
    limiting absorption principle for Schr\"odinger operators with Von
    Neumann-Wigner potentials.} Int. Eq. and Op. Theory, vol.\  14 (1991). 
%
\bibitem[FMS1]{fms1} J.\ Faupin, J.S.\ M\o ller, E.\ Skibsted: {\em
    Regularity of bound states.} arXiv:1006.5871. 
%
\bibitem[FMS2]{fms2} J.\ Faupin, J.S.\ M\o ller, E.\ Skibsted: {\em
    Second order perturbation theory for embedded eigenvalues.}
arXiv:1006.5869.
%
\bibitem[FH]{fh} R.G.\ Froese, I.\ Herbst: {\em Exponential bounds and   
absence of positive eigenvalues for $N$-body Schr\"odinger operators.}   
Comm.\ Math.\ Phys.\ 87, 429-447 (1982).   
%
\bibitem[FHHH1]{fhhh1} R.G.\ Froese, I.\ Herbst, M.\ Hoffmann-Ostenhof, T.\ 
Hoffmann-Ostenhof: {\em $L^2$-Exponential Lower Bounds to solutions
  of the Schr\"odinger Equation.}   
Comm.\ Math.\ Phys.\ 87, 265-286 (1982). 
%
\bibitem[FHHH2]{fhhh2} R.G.\ Froese, I.\ Herbst, M.\ Hoffmann-Ostenhof, T.\ 
Hoffmann-Ostenhof : {\em On the absence of
  positive eigenvalues for one-body Schr\"odinger operators.}   
J. Analyse Math.\ 41, 272-284, (1982). 
%
\bibitem[GG\'e]{gg} V.\ Georgescu, C.\ G\'erard: {\em On  the Virial
  Theorem in Quantum Mechanics}, Comm.\ Math.\ Phys.\ {\bf 208}, 275--281,
  (1999).  
%
\bibitem[GGM1]{ggm1} V.\ Georgescu, C.\ G\'erard, J.S.\ M\o ller: {\em Commutators, 
$C\sb 0$-semigroups and resolvent estimates}, J. Funct. Anal. {\bf
216}, no 2, p.\  303-361, 2004. 
%
%
\bibitem[GGo]{ggo} V.\ Georgescu, S.\ Gol\'enia: {\em Isometries, Fock spaces and spectral
analysis of Schr\"odinger operators on trees.}, J. Funct. Anal. 227 (2005), 389-429.
%
\bibitem[G\'e]{ge} C.\ G\'erard: {\em A proof of the abstract limiting
absorption principle by energy estimates}, J. Funct. Anal. 254 (2008), no 11, 2707-2724. 
%
\bibitem[GJ]{gj} S.\ Gol\'enia, Th.\ Jecko: {\em A new look at
  Mourre's commutator theory}, Compl. Ana. Op. Theory, Vol. 1, No. 3, 
p.\ 399-422, August 2007. 
%
\bibitem[GM]{GoleniaMoroianu}
S.\ Gol{\'e}nia, S.\ Moroianu: \emph{Spectral analysis of magnetic Laplacians
 on conformally cusp manifolds.}, Ann. H. Poincar{\'e} \textbf{9} (2008),
 no.~1, 131--179.
%
\bibitem[HeS]{hs} B.\ Helffer, J.\ Sj\"{o}strand: {\em Op\'erateurs de  
Schr\"{o}dinger avec champs magn\'etiques faibles et constants.}   
Expos\'e No.\ XII, S\'eminaire EDP, f\'evrier 1989, Ecole Polytechnique.  
%
\bibitem[H\"o1]{h1} L.\ H\"ormander: {\em The analysis of linear partial differential operators I.}, 
Springer-Verlag Berlin Heidelberg 1983. 
%
\bibitem[H\"o2]{h2} L.\ H\"ormander: {\em The analysis of linear partial differential operators II.}, 
Springer-Verlag Berlin Heidelberg 1983. 
%
\bibitem[H\"o3]{h3} L.\ H\"ormander: {\em The analysis of linear partial differential operators III.}, 
Springer-Verlag Berlin Heidelberg 1983. 
%
\bibitem[HuS]{hus} W.\ Hunziker, I.M.\ Sigal: {\em The quantum $N$-body
  problem}, J.\ Math.\ Phys.\ {\bf 41} (6), 3448--3510, 2000. 
%
\bibitem[IJ]{ij}A.D.\ Ionescu, D.\ Jerison: {\em On the absence of
  positive eigenvalues of Schr\"odinger operators with rough potentials},   
GAFA, Geom.\ Funct.\ Anal. 13, 2003, 1029-1081.  
%
\bibitem[Je1]{jec1}T.\ Jecko: {\em From classical to semiclassical non-trapping behaviour},   
C.\ R.\ Acad.\ Sci.\ Paris, Ser.\ I, {\bf 338}, p.\ 545--548, 2004.   
%
\bibitem[Je2]{jec2}T.\ Jecko: {\em Non-trapping condition for semiclassical 
Schr\"odinger operators with matrix-valued potentials.} Math. Phys. Electronic Journal, 
No. 2, vol. 11, 2005.
%
\bibitem[JMP]{jmp}A.\ Jensen, E.\ Mourre, P.\ Perry: {\em Multiple
  commutator  estimates and resolvent smoothness in quantum scattering
 theory.} Ann.\  Inst.\ H. Poincar\'e vol.\ 41, no 2, 1984, p.\
  207-225. 
%
\bibitem[Ka]{k} T. Kato: {\em Pertubation theory for linear operators.}
Springer-Verlag 1995.
%
\bibitem[Ki]{K} A.\ Kiselev: \emph{Imbedded singular continuous
    spectrum for Schr\"odinger operators},  J.\ Amer.\ Math.\ Soc.\ 18
  (2005), no.\ 3, 571–603.
%
\bibitem[Le]{l} N. Lerner: {\em Metrics on phase space and
    non-selfadjoint pseudodifferential  operators.} Birkha\"user 2010.
%
%
\bibitem[MU]{mu}K.\ Mochizuki, J.\ Uchiyama: {\em Radiation conditions
    and spectral  theory for 2-body Schr\"odinger operators with
    ``oscillating'' long-range potentials  I: the principle of
    limiting absorption.} J.\ Math.\ Kyoto Univ.\ {\bf 18}, 2,  
377--408, 1978. 
%
\bibitem[M\o]{mo} J. S. M\o ller: {\em An abstract radiaton condition and applications to
N-body systems.}, Reviews of Math. Phys., Vol. 12, No. 5 (2000), 767-803.

%
\bibitem[MS]{ms} J.S.\ M\o ller, E.\ Skibsted:
\emph{Spectral theory of time-periodic many-body systems.}
Adv.\ Math.\ {\bf 188} (2004), no.\ 1, 137--221.
%

\bibitem[Mo]{m}E.\ Mourre: {\em Absence of singular continuous
    spectrum for    
certain self-adjoint operators.} Comm.\ in Math.\ Phys.\ {\bf 78}, 391--408,
  1981.  
%
%
%
\bibitem[RS4]{rs4}M.\ Reed, B.\ Simon: {\em Methods of Modern
    Mathematical  Physics, Tome IV: Analysis of operators.} Academic Press.  
%
\bibitem[ReT1]{ret1}P.\ Rejto, M.\ Taboada: {\em A limiting absorption
    principle for Schr\"odinger operators with generalized Von
    Neumann-Wigner potentials I. Construction of approximate phase.}
  J.\ Math. Anal. and Appl. 208, p.\ 85--108 (1997). 
%
\bibitem[ReT2]{ret2}P.\ Rejto, M.\ Taboada: {\em A limiting
    absorption principle for Schr\"odinger operators with generalized
    Von Neumann-Wigner potentials II. The proof. } J. Math.
  Anal. and Appl. 208, p.\ 311-336 (1997).
%
\bibitem[Re]{r} C.\ Remling: \emph{The absolutely continuous spectrum
    of one-dimensional Schr\"odinger operators with decaying
    potentials}, Comm.\ Math.\ Phys.\ {\bf 193} (1998), 151--170.
%
\bibitem[RoT]{rt}D.\ Robert, H. Tamura: {\em Semiclassical estimates for 
resolvents and asymptotics for total cross-section.} Ann.\  Inst.\ H. Poincar\'e {\bf 46}, 
415-442 (1987).
%
\bibitem[Roy]{roy}J.\ Royer: {\em Limiting absorption principle for
    the dissipative Helmholtz  equation.} CPDE {\bf 35}, 1458-1489 (2010).
%
\bibitem[Sa]{s}J.\ Sahbani: {\em  The conjugate operator method for locally 
  regular Hamiltonians.} J.\ Oper.\ Theory 38, No.\ 2, 297--322
  (1997).
 %
\bibitem[Si]{si}B.\ Simon: {\em  On positive eigenvalues of one body
  Schr\"odinger operators.} Comm.\ Pure Appl.\ Math.\ 22, 531--538 (1969).
%
\bibitem[Wa]{w}X.P. Wang: {\em Semiclassical resolvent estimates for 
$N$-body Schr\"odinger operators.} J.\ Funct.\ Anal.\
{\bf  97}, 466--483 (1991).
%
\bibitem[WZ]{wz}X.P. Wang, P. Zhang: {\em High frequency limit of the 
Helmholtz equation with variable refraction index.} J.\ Funct.\ Anal.\
{\bf 230}, 116-168 (2006). 
%
\end{thebibliography}
\end{document}